\documentclass[11pt]{article}
\usepackage[margin=1.3in]{geometry}
\usepackage{graphicx} % Required for inserting images
\usepackage{authblk}  % For better author formatting
\usepackage{hyperref} % For clickable email links
\usepackage[dvipsnames]{xcolor}
\usepackage{esvect}
\usepackage[ruled,vlined,linesnumbered]{algorithm2e}
\SetAlgoNoLine
\usepackage{graphicx}      
\usepackage{amsmath,amssymb, amsfonts,amsthm}

\usepackage{subcaption}
\usepackage{standalone}
\usepackage{hyperref}   
\usepackage{color}% for hyperlinks
\newtheorem{Theorem}{Theorem}
\newtheorem{Lemma}{Lemma}
\newtheorem{Remark}{Remark}
\newtheorem{proposition}[Theorem]{Proposition} 
\usepackage{tikz}
\usetikzlibrary{arrows.meta,calc}
\usepackage{float}
\usepackage{todonotes}
\usepackage[dvipsnames]{xcolor}
\usepackage{enumitem}
\usepackage{url}
\usepackage{array}
\usepackage{caption}
\usepackage{amsmath}
\usepackage{nicematrix}
\usepackage{arydshln}  
\usepackage{array} % Required for tabular
\usepackage{amsfonts}
\usepackage{amsmath,amssymb}
\usepackage{amsthm}
\usepackage{graphicx}      % include this line if your document contains figures
\usepackage{comment}

\usepackage{tikz}
\usepackage{pgfplots}
\pgfplotsset{compat=1.8}
\usepackage{pgfplotstable}
\usetikzlibrary{3d, calc}
\usepackage{placeins}   % Ensures figures do not float past barriers
\usepackage{cases}
\pgfplotsset{compat=1.17}

\usepackage[
    backend = biber,
    doi = false,
    url = false,
    isbn = false,
    sorting = none,
    giveninits = true,  % modern replacement for firstinits
    maxnames = 99,
    hyperref = true,
    block = none
]{biblatex}

\renewbibmacro{in:}{}

\bibliography{referenece}

\newtheorem{Definition}{Definition}

\newcommand{\N}{\mathbb{N}}
\newcommand{\iN}{N_{\rm in}}
\newcommand{\bS}{\mathbb{S}}

\newcommand{\cE}{\mathcal{E}}
\newcommand{\R}{\mathbb{R}}

\newcommand{\cV}{\mathcal{V}}
\newcommand{\cG}{\mathcal{G}}
\newcommand{\cF}{\mathcal{F}}
\newcommand{\xc}[1]{\vspace{.1cm}

\noindent {\em #1} }

\makeatletter
\def\blfootnote{\xdef\@thefnmark{}\@footnotetext}
\makeatother

\title{Switched Linear Ensemble Systems and \\ Structural Controllability}

\begin{document}
\author{ Haoyu Yin, \quad Yi Li, \quad Ouyang Du, \quad Bruno Sinopoli, \quad and \quad Xudong Chen$^*$ 
}
\date{}
\maketitle
\blfootnote{All the authors are with the Electrical and Systems Engineering, Washington University in St. Louis. Emails: \texttt{\{h.yin, lyi1, d.ouyang, bsinopoli, cxudong\}@wustl.edu}. }

\blfootnote{$^*$Corresponding author: X. Chen.}
%\blfootnote{This work was supported by NSF ECCS-2426017 and by DARPA HR0011-25-3-0211.}

\begin{abstract}
This paper introduces and solves a structural controllability problem for ensembles of switched linear systems. All individual systems in the ensemble are sparse and governed by the same sparsity pattern, and undergo switching among subsystems by following the same switching sequence. The controllability of an ensemble system describes the ability to use a common control input to simultaneously steer every individual system. A sparsity pattern is called structurally controllable for pair \((k,q)\) if it admits a controllable ensemble of \(q\) individual systems with at most \(k\) subsystems. We derive a necessary and sufficient condition for a sparsity pattern to be structurally controllable for a given \((k,q)\), and characterize when a sparsity pattern admits a finite \(k\) that guarantees structural controllability for \((k,q)\) for arbitrary $q$.
 Compared with the linear time-invariant ensemble case, this second condition is strictly weaker. We further show that these conditions have natural connections with maximum flow, and hence can be checked by polynomial algorithms. Specifically, the time complexity of deciding structural controllability is \(O(n^3)\) and the complexity of computing { the smallest number of subsystems needed} is \(O(n^3 \log n)\), with \(n\) the dimension of each individual system. %The results we present in this paper provide scalable tools and insights for controlling large ensembles of switched systems.
\end{abstract}

% \begin{keywords}
% Structural system theory, linear ensemble systems, switched systems, graph theory, maximum flow.
% \end{keywords}
\section{Introduction}
An ensemble system comprises a large number of individual systems, and ensemble control deals with the problem of using a common, finite-dimensional control input to simultaneously steer all the individual systems. 
It originated from quantum control~\cite{brockett2000stochastic,li2006control,chen2020ensemble} where one uses radio frequency impulses to steer an ensemble of spin dynamics, and finds applications in many other disciplines, ranging from neuroscience~\cite{ching2013control,zlotnik2016phase}, learning and neural ODEs~\cite{tabuada2022universal,bayram2025geometric}, to multi-agent control~\cite{becker2017controlling,chen2019controllability}. 

A distinctive feature of ensemble control is that the dynamics of the individual systems are independent of each other (precisely, the dynamics of any given individual system do not depend on the states of the others),  
and the coupling is through the common control input only. 
Controllability of the whole ensemble is attained by variations in parameters of the individual systems~\cite{li2006control,helmke2014uniform,chen2019structure}.  
This is in contrast with network control where controllability of a networked system relies often on a certain leader-follower hierarchy. 
Ensemble control thus provides an alternative approach for controlling large-scale systems: Instead of controlling a large networked system, one controls a population of small ones, which is by nature resilient and scalable; indeed, malfunction of a few individual systems do not affect the others and one is free to add into or remove from the ensemble an arbitrary number of individual systems~\cite{chen2023sparse}.

A main challenge of ensemble control stems from the requirement that the control input be generated irrespective of the individual system. Controllability of each individual system is necessary for the ensemble to be controllable, yet far from being sufficient. For example, an ensemble of integrators is not controllable as a whole, regardless of their respective actuation gains. Necessary and/or sufficient conditions have been obtained for ensemble controllability, yet open questions still remain. See, for example,~\cite{dirr2021uniform,chen2023controllability,chen2019structure} and the references therein.

In the recent work~\cite{chen2023sparse}, we considered continuum ensembles of linear {\it time-invariant} (LTI) systems, whose $(A, B)$ pairs share a common sparsity pattern.  
We investigated what type of sparsity pattern can sustain controllability of such ensemble system, i.e., we address the structural controllability problem for linear ensemble systems.    
Amongst others, we have provided a necessary and sufficient condition~\cite[Theorem 1]{chen2023sparse} for the sparsity pattern to be structurally ensemble controllable. Without recalling the full and precise statement, we note here that if a sparsity pattern is structurally ensemble controllable, then the digraph corresponding to the sparsity pattern of the $A$-matrix must have a node-wise disjoint cycle cover (node-wise Hamiltonian decomposition). 
In particular, the condition for structural ensemble controllability is strictly stronger than the one for structural controllability (for LTI systems).  

In this paper, we follow the research line of~\cite{chen2023sparse} and continue to investigate ensemble controllability from the perspective of structural system theory. 
A major motivation of the current work is rooted in the following question: To what extent can we relax the condition established in~\cite{chen2023sparse}? 
Said in other words, if a sparsity pattern $\bS$ cannot sustain controllability for any LTI ensemble system, then is it possible that we find another class of ensemble systems that are controllable and, at the same time, compliant with the sparsity pattern? 
We provide an affirmative answer to the question.  

Specifically, we consider in this paper the class of {\it switched}, linear ensemble systems and address the associated structural controllability problem. The main contributions of this work are the following:    
\begin{enumerate}
\item Given a pair $(k,q)$ of positive integers, where {  $k$ is the number of subsystems} and $q$ is the number of individual systems, we provide in Theorem~\ref{thm:main} a necessary and sufficient condition for a sparsity pattern $\bS$ to sustain controllability for switched, linear $q$-ensemble systems (i.e., an ensemble of $q$ switched, linear systems sharing a common switching signal). 
For convenience, we call any such sparsity pattern $\bS$ structurally controllable for $(k,q)$.  
\item Leveraging the above result, we provide in Proposition~\ref{prop:KS} a necessary and sufficient condition for a sparsity pattern $\bS$ to admit a finite $k$ such that $\bS$ is structurally controllable for $(k,q)$ for any $q\in \N$ (with $k$ fixed). Compared with the condition for LTI ensemble systems, the condition for the switched, linear ensemble systems is strictly weaker, requiring only that sparsity pattern of the $A$-matrix be such that every column of $A$ contains at least one $\star$-entry. The same condition, said in graphical terms, is that the digraph corresponding to the $A$-matrix has the property that each node has at least one in-neighbor.      
\item In addition to the two necessary and sufficient conditions, we further present in Section~\ref{ssec:algorithms} polynomial algorithms for checking the conditions. This is done by relating the conditions to appropriately defined maximum flow problems (Theorem~\ref{thm:main2}).  
The time complexity of checking whether a sparsity pattern is structurally controllable for a given pair $(k,q)$ is $O(n^3)$ where $n$ is the dimension of the {\it individual} system (not the dimension of the ensemble system), independent of $k$ or $q$ (note that the setting of our interest is such that $q \gg 1$). 
The time complexity for checking whether there exists a finite $k$ such that $\bS$ is structurally controllable for $(k,q)$ for any $q\in \N$ is also $O(n^3)$. Finally, if such $k$ exists, then the time complexity for computing the smallest value $k^*$ is $O(n^3\log n)$.    
\end{enumerate}

Since the seminal work of Lin~\cite{lin1974structural}, there have been many efforts and progress made in structural system theory. 
To the best of our knowledge, there has not been any work in the literature that addresses structural controllability for switched, ensemble systems.  
We mention below a few works that are relevant to this paper. We will focus on two areas, namely, structural system theory in ensemble systems and in switched systems. The problem of structural ensemble controllability was initiated and addressed in~\cite{chen2023sparse} as mentioned above. More recently, the problem of structural averaged controllability has been addressed in~\cite{gharesifard2021structural,chen2025structural}, where the goal of that research line is to understand to what extent can we relax the condition on the sparsity pattern if only certain output functions (e.g., average and higher moments of an ensemble) are of one's interest? The class of ensemble dynamics considered in the above work are LTI ensemble systems. 
Next, for the relevant work about (finite dimensional) switched systems, we first mention~\cite{liu2013structural} where the authors considered switched linear systems, with heterogeneous sparsity patterns over different time periods. Given the number~$k$ of subsystems, the authors provided a necessary and sufficient condition (see Theorems~9 and~11 of their work) for the sequence of sparsity patterns to be structurally controllable. {  A more rigorous and complete proof of the same result was later established in~\cite{zhang2024generalized}.}
The authors of~\cite{liu2013structural} have also articulated an algorithm for checking the condition, with time complexity $O(\sqrt{kn} kn^2)$.    
We also mention~\cite{Cui2020,cui2022temporal} where the authors considered a similar setting, but with the focus on computing the smallest value $k^*$ (note that such $k^*$ always exists for a single linear system, but is not the case for the ensemble setting). The authors have also provided an algorithm to obtain $k^*$ with time complexity $O(n^4)$. 
{  
Finally, we note that there has been a different structural controllability problem formulated for switched, linear system with fixed switching sequence~\cite{zhang2024reachability, lebon2024controllability}, in which 
the authors fix the number of allowable switches, and treat the state matrices over different time segments (between switches) as free parameters constrained by given sparsity patterns.}

The remainder of this paper is organized as follows. In Section~\ref{sec:notation}, we gather key notions and preliminaries for our discussion.
The main results of the paper are presented in Section~\ref{sec: Problem and Result}, in which we present the two necessary and sufficient conditions (Theorem~\ref{thm:main} and Proposition~\ref{prop:KS}), connections with maximum flow (Theorem~\ref{thm:main2}), and the polynomial algorithms to check the conditions. We provide the proof of Theorem~\ref{thm:main} in Section~\ref{sec:proof_main}, and the proof of Theorem~\ref{thm:main2} in Section~\ref{sec:proof_flow_equiv}. The paper ends with conclusions and outlooks.

\section{Preliminaries}\label{sec:notation}

\subsection{Switched linear ensemble systems}
{ 
We consider in this paper finite ensembles of switched, linear systems:
\begin{equation}\label{eq:feltv}
\dot x_p(t) = A_p[\sigma(t)] x_p(t) + B_p[\sigma(t)] u(t), \mbox{ for } p = 1,\ldots,q,
\end{equation}
where $x_p(t)\in \R^n$ is the state of the $p$th individual system, $u(t)\in \R^m$ is the common control input, $\sigma: [0,\infty)\to \{1,\ldots, k\}$ is the common switching signal (which is a right-continuous, piecewise constant function), and $$\begin{bmatrix} A_p[\ell] & B_p[\ell] \end{bmatrix}\in \R^{n\times (n + m)} \quad \mbox{for all }\ell = 1,\ldots, k.$$  
We call~\eqref{eq:feltv} a {\it switched, linear $q$-ensemble system with $k$ subsystems}. 
For convenience, let $x(t):= (x_1(t);\ldots;x_q(t))$ and 
$$
A[\ell]:=  
\begin{bmatrix}
A_1[\ell] & & \\
& \ddots & \\
& & A_q[\ell]
\end{bmatrix} 
\quad \mbox{and} \quad 
B[\ell]:=
\begin{bmatrix}
B_1[\ell] \\
\vdots \\
B_q[\ell]
\end{bmatrix}.
$$
Then, system~\eqref{eq:feltv} can be written as
\begin{equation}
\label{eq:kswitchsys}
\dot x(t) = A[\sigma(t)] x(t) + B[\sigma(t)] u(t).
\end{equation}

The switched system~\eqref{eq:kswitchsys} is said to be controllable over the time interval $[0,T]$ if for any
initial state $x(0)\in\mathbb{R}^{qn}$ and any target state
$x^*\in\mathbb{R}^{qn}$, there exist a switching signal $\sigma:[0,T]\to\{1,\dots,k\}$ and an integrable control input
$u:[0,T]\to\mathbb{R}^m$ such that the solution $x(t)$ generated by~\eqref{eq:kswitchsys} satisfies $x(T) = x^*$.
% \end{Definition}
}

% Let 
% $$x(t):= 
% \begin{bmatrix}
% x_1(t) \\
% \vdots \\
% x_q(t)
% \end{bmatrix}, \quad 
% A(t):= 
% \begin{bmatrix}
% A_1(t) & & \\
% & \ddots & \\
% & & A_q(t)
% \end{bmatrix}, \quad \mbox{and} \quad
% B(t):=
% \begin{bmatrix}
% B_1(t) \\
% \vdots \\
% B_q(t)
% \end{bmatrix}, 
% $$
% where $A(t)$ is a block-diagonal matrix. Then, 
% one can re-write~\eqref{eq:feltv} as a $qn$-dimensional linear time-varying control system: 
% \begin{equation}\label{eq:ltvsys}
% \dot x(t) = A(t) x(t) + B(t) u(t).
% \end{equation}
% We say that system~\eqref{eq:ltvsys} is {\it controllable over $[0,T]$}, for $T > 0$, if for any pair of initial state $x(0)$ and target state $x^*$ at time~$T$, there exists an integrable function $u:[0,T]\to \R^m$ such that the solution of~\eqref{eq:ltvsys} generated by $u(t)$ satisfies $x(T) = x^*$.   

% If the matrix pair $(A(t), B(t))$ is piecewise constant in $t$, then the system~\eqref{eq:ltvsys} is a switched linear system.
{  Controllability of the switched linear system~\eqref{eq:kswitchsys} has been widely addressed in the literature (see, e.g.,~\cite{SUN2002,liberzon2003switching,liu2013structural}). We recall the following result adapted from~\cite{SUN2002}:

\begin{Lemma}
\label{lem:ctrlkswitchsys}
Let $W$ be the smallest subspace of $\R^{nq}$ that contains the column space of each $B[\ell]$ and is invariant under $A[\ell]$ (i.e., for any $w\in W$, $A[\ell]w \in W$) for all $\ell = 1,\ldots, k$. Then, the switched  system~\eqref{eq:kswitchsys}, with  $(A[\ell],  B[\ell])$'s the matrix pairs for the $k$ subsystems, is controllable over $[0,T]$, for some (and hence, any) $T > 0$, if and only if $W = \R^{nq}$.
\end{Lemma}

%     $$\left\langle
% A[1],\dots,A[k]\mid B[1],\dots,B[k]
%     \right\rangle
%     =
%     \mathbb R^{qn},
%     $$
%     where
%     the left hand side of the equation
%     denotes the smallest subspace of $\mathbb R^{qn}$ containing
%     $\sum_{\ell=1}^{k}\operatorname{im} B[\ell]$ and invariant under
%     $A[1],\dots,A[k]$.

% Under the definition, one does not fix a switching
% time sequence in advance. Instead, controllability is defined by the existence of a
% finite terminal time, a piecewise-constant switching signal, and a control input
% that steer the system from a given initial state to a prescribed target state.
% Accordingly, the controllability criterion in Lemma~\ref{lem:ctrlkswitchsys}
% depends only on the subsystem family $\{(A[\ell],B[\ell])\}_{\ell=1}^{k}$ and not
% on any given switching instants.
}

\subsection{Sparsity pattern, graph representation, and maximum flow}
In this paper, a {\it sparsity pattern} $\bS$ refers to a vector subspace of $\R^{n\times (n + m)}$ such that it admits a basis consisting only of matrices $E_{ij}$'s, i.e., matrices with $1$ on the $(i,j)$th entry and $0$ elsewhere. 
One can represent $\bS$ by a matrix with $\star$- and $0$-entries.  
The $(i,j)$th entry is a $\star$-entry if $E_{ij}\in \bS$ and is, otherwise, a $0$-entry.

The sparsity pattern $\bS\subseteq\R^{n\times (n + m)}$ induces a directed graph $G=(V,E)$ whose node and edge sets are defined as follows:
\begin{enumerate}
\item The node set $V$ is a union of two sets  $V^\alpha := \{\alpha_1,\ldots, \alpha_n\}$ and $V^\beta := \{\beta_1,\ldots, \beta_m\}$. We call the nodes in $V^\alpha$ the {\it state nodes} and the nodes in $V^\beta$ the {\it control nodes};
\item The edge set $E$ is defined as follows: There exists an edge from $\alpha_j$ to $\alpha_i$ (resp., from $\beta_j$ to $\alpha_i$) if the $(i,j)$th (resp., $(i,j+n)$th) entry of $\bS$ is a $\star$-entry.   
\end{enumerate}

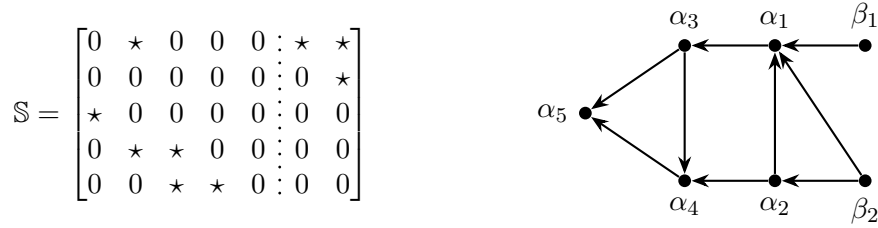
\begin{figure}[ht]
    % \centering
    \begin{minipage}{0.45\linewidth}
        \[ \mathbb{S}=
        \left[
\begin{NiceArray}{ccccc:cc}
  0 & \star & 0 & 0 & 0 & \star & \star\\
  0 & 0     & 0 & 0 & 0 & 0     & \star\\
  \star & 0 & 0 & 0 & 0 & 0     & 0\\
  0 & \star & \star & 0 & 0 & 0 & 0\\
  0 & 0     & \star & \star & 0 & 0 & 0
\end{NiceArray}
        \right]
        \]
    \end{minipage} \quad \quad 
    \begin{minipage}{0.5\linewidth}
        % \centering
        \begin{tikzpicture}[>=Stealth, scale=0.98, rotate=-90]
  \tikzset{dot/.style={circle,fill,inner sep=1.7pt}};
  %--- nodes
  \node[dot,label=below:$\beta_2$]  (beta_2) at (2,2)   {};
  \node[dot,label=above:$\alpha_1$]  (alpha_1) at (0.5,1)  {};
  \node[dot,label=below:$\alpha_2$]   (alpha_2) at (2,1)  {};

  \node[dot,label=above:$\alpha_3$]   (alpha_3) at (0.5,0)  {};
  \node[left of=alpha_3,xshift=-0.75cm, yshift=0.25cm](name){};\node[dot,label=above:$\beta_1$]  (beta_1) at (0.5,2)   {};
  \node[dot,label=below:$\alpha_4$]   (alpha_4) at (2,0)  {};

  \node[dot,label=left:$\alpha_5$]   (alpha_5) at (1.25,-1.1)  {};
  %--- curved bidirectional arrows
  \draw[->,thick] (beta_1) to(alpha_1);
  \draw[->,thick] (beta_2) to (alpha_2);
  \draw[->,thick] (beta_2) to (alpha_1);
  \draw[->,thick] (alpha_2) -- (alpha_1);
  \draw[->,thick] (alpha_1) to (alpha_3);
  \draw[->,thick] (alpha_2) to (alpha_4);
  \draw[->,thick] (alpha_3) to (alpha_5);
  \draw[->,thick] (alpha_4) to (alpha_5);
  \draw[->,thick] (alpha_3) -- (alpha_4);
  % \draw[->,thick] (alpha_3) -- (alpha_4);
  %--- self loop at u4
\end{tikzpicture}
    \end{minipage}
    \caption{{\it Left:} A sparsity pattern $\bS$. {\it Right:} The digraph $G$ associated with $\bS$.}
    \label{Fig1}
\end{figure}

% its associated digraph is defined as
% \[
%     G = (V, E),
% \]
% where the node set $V$ is the disjoint union of the state-node set
% \[
%     V^{\alpha} := \{\alpha_1,\ldots,\alpha_n\},
% \]
% and the input-node set
% \[
%     V^{\beta} := \{\beta_1,\ldots,\beta_m\}.
% \]
% The edge set $E$ is given by
% \[
%     E := \{\, (\alpha_i,\alpha_j) \mid A_{ji} = \star \,\}
%      \cup 
%     \{\, (\beta_i,\alpha_j) \mid B_{ji} = \star \,\}.
% \]
% In particular, there is a directed edge from $\alpha_i$ to $\alpha_j$ whenever the entry in row $j$, column $i$ of $A$ is structurally nonzero, and a directed edge from $\beta_i$ to $\alpha_j$ whenever the entry in row $j$, column $i$ of $B$ is structurally nonzero. 
An example of a sparsity pattern and its associated digraph for $n=5$ and $m=2$ is shown in Fig.~\ref{Fig1}.

%We then introduce the in- and out-neighbor sets associated with the digraph $G$.

%Let $G = (V,E)$ be the digraph associated with the sparsity pattern $\mathbb{S}$, with
%$V = V^{\alpha} \cup V^{\beta}$ as defined above. 
For any subset $V' \subseteq V^{\alpha}$, the {\it in-neighbor set} of $V'$ in $G$ is defined as
\[
    N_{\mathrm{in}}(V')
    := \bigl\{\, v \in V \,\bigm|\, \exists\, w \in V' \text{ such that } (v,w) \in E \bigr\}.
\]
We further define  
\[
    N_{\mathrm{in}}^{\alpha}(V') := N_{\mathrm{in}}(V') \cap V^{\alpha}
    \quad \mbox{and} \quad
    N_{\mathrm{in}}^{\beta}(V') := N_{\mathrm{in}}(V') \cap V^{\beta}.
\]
Similarly, the {\it out-neighbor set} of $V'$ in $G$ is defined as
\[
    N_{\mathrm{out}}(V')
    := \bigl\{\, w \in V \,\bigm|\, \exists\, v \in V' \text{ such that } (v,w) \in E \bigr\}.
\]

% Recall that a {\bf bipartite graph} is a graph $G = (V_L, V_R, E_{LR})$ whose vertex set is partitioned into two
% disjoint subsets $V_L$ and $V_R$ such that every edge in $E_{LR}$ connects a vertex in $V_L$ to a
% vertex in $V_R$. No edge is allowed between two vertices in the same subset.

% Given any digraph $G = (V,E)$ with input-node set $V^{\beta} = \{\beta_1,\ldots,\beta_m\}$ and
% state-node set $V^{\alpha} = \{\alpha_1,\ldots,\alpha_n\}$, 
We now introduce a class of maximum flow problems that are relevant to the main results. 
Given the digraph $G$ associated with $\bS$, we construct a weighted digraph $\cG = (\cV, \cE, c)$ on which a maximum flow problem is defined.
Let
\[
    \mathcal{V} := \{s,t\} \cup V_L \cup V_R,
\]
where $s$ and $t$ are the source and sink nodes, respectively, and
\[
    V_L := \{\lambda_1,\ldots,\lambda_m, \nu_1,\ldots,\nu_n\},
    \qquad
    V_R := \{\mu_1,\ldots,\mu_n\}.
\]
Here, $\lambda_i$ is a copy of the control node $\beta_i$, while $\nu_j$ and $\mu_j$ are two disjoint copies
of the state node $\alpha_j$, for $i\in\{1,\ldots,m\}$ and $j\in\{1,\ldots,n\}$.
We define the edge set as
\[
    \mathcal{E} := E_{sL} \cup E_{LR} \cup E_{Rt},
\]
where
\begin{align*}
    E_{LR}
    &:= \bigl\{ (\lambda_i,\mu_j) | \, (\beta_i,\alpha_j) \in E \bigr\}
    \cup
       \bigl\{(\nu_i,\mu_j) | \, (\alpha_i,\alpha_j) \in E \bigr\},\\
    E_{sL} &:= \{\, (s,v)\mid v\in V_L \,\},\\
    E_{Rt} &:= \{\, (v,t)\mid v\in V_R \,\}.
\end{align*}

Each edge $e\in\mathcal{E}$ is assigned a capacity $c(e)$, the value of which will be specified later.

Given the weighted digraph $\cG = (\cV,\cE,c)$, a function
$f:\cE \to \R_{\geq 0}$ is a {\it flow} on $\cG$ if it satisfies:
\begin{enumerate}
\item {\it Capacity constraint:} For every edge $e \in \cE$,
\begin{equation}
\label{eq:capacity_con}
    f(e) \leq c(e).
\end{equation}
\item {\it Conservation of flow:} For every node $v \in \cV \setminus \{s,t\}$,
\begin{equation}
\label{eq:flow_balance}
    \sum_{u \in N_{\rm in}(v)} f(u,v)
    = \sum_{w \in N_{\rm out}(v)} f(v,w).
\end{equation}
\end{enumerate}
The \emph{value} of a flow $f$ on $\cG$ is defined as
\begin{equation}\label{eq:value}
    |f| := \sum_{v \in N_{\rm out}(s)} f(s,v)= \sum_{w \in N_{\rm in}(t)} f(w,t).
\end{equation}
Let $\cF$ denote the set of all flows on $\cG$ satisfying
\eqref{eq:capacity_con} and \eqref{eq:flow_balance}. The so-called maximum flow problem~\cite{cormen2022introduction} is the following optimization problem:
\[
    \max_{f \in \cF} |f|,
\]
which can be solved in polynomial time.
% By~\cite{cormen2022introduction}, if there exists a flow $f^\star \in \cF$ such that $|f^\star(\cG)| = |V_R|$, then the bipartite graph
% $(V_L, V_R, E_{LR})$ admits a right perfect matching.

\section{Main Results} \label{sec: Problem and Result}
In this section, we present the main results. 
In Subsection~\ref{ssec:nec_suf_cond}, we formulate the structural controllability problem for switched, linear ensemble systems and provide in 
Theorem~\ref{thm:main} a necessary and sufficient condition for a sparsity pattern $\bS$ to be structurally controllable for a given pair $(k, q)$, with $k$ the number of subsystems and $q$ the number of individual systems. 
In Subsection~\ref{ssec:KS}, we present in Proposition~\ref{prop:KS} a necessary and sufficient condition for the existence of a finite $k$ such that $\bS$ is structurally controllable for $(k, q)$ for all $q \in \mathbb{N}$. 
Then, in Subsection~\ref{ssec:algorithms}, we relate the two conditions to some maximum flow problem, which enables us to use polynomial algorithms to check the two conditions and to compute the smallest value~$k^*$. 

\subsection{Necessary and sufficient condition}
We start with the following definition: 
\label{ssec:nec_suf_cond}

\begin{Definition}\label{def:structrl}
    A sparsity pattern $\bS\subseteq \R^{n\times (n + m)}$ is {\it structurally controllable} for $(k,q)$, with $k \geq 1$ and $q\geq 1$, if there exist matrices $\begin{bmatrix}A_p[\ell] &  B_p[\ell]\end{bmatrix}\in \bS$,  for all  $\ell = 1,\ldots, k$ and for all $p = 1,\ldots, q$ 
    such that system~\eqref{eq:kswitchsys} is controllable. 
\end{Definition}

Note that if $G$ is structurally controllable for some $(k,q)$, then by the above definition and Lemma~\ref{lem:ctrlkswitchsys} it is structurally controllable for any $(k',q')$ for $k' \geq k$ and $q'\leq q$.

The first main result of the paper provides a necessary and sufficient condition  for a sparsity pattern to be structurally controllable. 

\begin{Theorem}\label{thm:main}
Let $\bS\subseteq \R^{n\times (n + m)}$ be a sparsity pattern and $G$ be the associated digraph. Then, $G$ is structurally controllable for a given $(k,q) \in \mathbb{N} \times \mathbb{N}$ if and only if the following two items hold:
\begin{enumerate}
    \item For every state node $\alpha_i$, there exist a control node $\beta_j$ and a path from $\beta_j$ to $\alpha_i$;
    \item For any subset $V'\subseteq V^\alpha$,
    \vspace{-3pt}
\begin{equation}\label{eq:core}
    k |\iN^\beta(V')| + kq |\iN^\alpha(V')| \geq q |V'|. 
    \end{equation}

\end{enumerate}
\end{Theorem}

The above result generalizes known results in the literature for the special case where $q = 1$. First, note that if $(k,q) = (1,1)$, then the problem of structural controllability is reduced to the standard one for linear time-invariant systems~\cite{lin1974structural}. In particular, item~2 of Theorem~\ref{thm:main} is nothing but  
$|\iN(V')| \geq |V'|$ for all $V'\subseteq V^\alpha$~\cite{olshevsky2015minimum}. 
Next, note that for $q = 1$ but with $k$ arbitrary, the structural controllability problem has been investigated in~\cite{liu2013structural}, where the authors derived a necessary and sufficient condition (see, Theorem~9 of their paper) using the so-called colored union graphs {  (see, also,~\cite{zhang2024generalized} for a complete proof of the same result).}  

%We show in Appendix~\ref{sec:appendix-proof-lemma3} that Theorem~\ref{thm:main} is equivalent with their result.  

Next, we specialize Theorem~\ref{thm:main} to the case $k = 1$, i.e., we deal with structural controllability problem for time-invariant linear ensemble systems. The second item of { Theorem~\ref{thm:main}} is then reduced to
\begin{equation}\label{eq:k=1}
 |\iN^\beta(V')| + q |\iN^\alpha(V')| \geq q |V'|, \quad \mbox{for all } V'\subseteq V^\alpha.
\end{equation}
Note, in particular, that~\eqref{eq:k=1} holds for {\it any} sufficiently large $q$ if and only if $$|\iN^\alpha(V')| \geq |V'|, \quad \mbox{for all } V'\subseteq V^\alpha.$$ 
It is well known (as a consequence of the Hall's marriage theorem) that the above condition holds if and only if the subgraph $G^\alpha$ induced by $V^\alpha$ has a node-wise disjoint cycle cover. This condition, together with item~1 of Theorem~\ref{thm:main}, coincides with the necessary and sufficient condition obtained in~\cite{chen2023sparse} which deals with the structural controllability problem for continuum of linear ensemble systems.

\subsection{Toward infinite $q$ with finite $k$}
\label{ssec:KS}
{ 
In this subsection, we address the following question: {\it When can a sparsity pattern $\mathbb{S}$ sustain an arbitrarily large number $q$ of individual systems with only a fixed number $k$ of subsystems?} 
}
Specifically, for a given $\mathbb{S}$, we define
\begin{equation}
\label{eq:KS}
\mathcal{K}(\mathbb{S}):=\{k\in \N \mid G \mbox{ is structurally controllable }  \text{for } (k,q) \mbox{ for any } q\in \N \}. 
\end{equation}
We provide below a complete characterization of the set $\mathcal{K}(\mathbb{S})$, which follows as a consequence of Theorem~\ref{thm:main}:

\begin{proposition}\label{prop:KS}
Given a sparsity pattern  $\mathbb{S}\subseteq \R^{n\times (n + m)}$, 
let $G$ be the associated digraph and $\mathcal{K}(\mathbb{S})$ be given as in~\eqref{eq:KS}.  
Let $k^*\in \N\cup \{\infty\}$ be given as follows:
\begin{equation}
\label{eq:k_star}
k^*:= \max_{V' \subseteq V^\alpha} \left\lceil \frac{|V'|}{N^\alpha_{\mathrm{in}}(V')} \right\rceil.
\end{equation}
Then, $\mathcal{K}(\mathbb{S}) = \{ k \in \N \mid k\geq k^*\}$. 
\end{proposition}

\begin{proof}
We re-write~\eqref{eq:core} as follows
\begin{equation}
\label{eq:cond_k_leq}
    k\geq \frac{|V'|}{|\iN^\beta(V')|/q+|\iN^\alpha(V')|},
\end{equation}
which must hold for all $q \in \mathbb{N}$,
where the right hand side is infinity if $N_{\rm in}(V')$ is empty. 
Note that the denominator is monotonically decreasing in $q$ and converges to $|N^\alpha_{\rm in}(V')|$ as $q\to\infty$. It follows that if~\eqref{eq:cond_k_leq} holds for all $q\in \N$, then we must have that 
$$k\geq  \frac{|V'|}{|\iN^\alpha(V')|}.$$
It follows that the minimal $k\in \N\cup\{\infty\}$ satisfying the above equation for any subset $V' \subseteq V^\alpha$ is $k^*$ given in~\eqref{eq:k_star}.
\end{proof}

\begin{Remark}\normalfont
By Proposition~\ref{prop:KS}, $\mathcal{K}(\mathbb{S})$ is nonempty if and only if $k^*$ is finite, which holds if and only if $N^\alpha_{\mathrm{in}}(V')\neq \varnothing$ for any $V'\subseteq V^\alpha$. 
Moreover, if $k^*$ is finite, then it is bounded above by~$n$.\hfill{\qed}
\end{Remark}

\subsection{Polynomial algorithms for Theorem~\ref{thm:main} and Proposition~\ref{prop:KS}}
\label{ssec:algorithms}

In this subsection, we present a polynomial algorithm, with time complexity $O(|V|^3)$, for checking the necessary and sufficient condition in Theorem~\ref{thm:main}, and a polynomial algorithm, with time complexity $O(|V|^3\log|V|)$, to compute the value of $k^*$ introduced in Proposition~\ref{prop:KS}.

\subsubsection{Algorithm for Theorem~\ref{thm:main}}
Note that the first item of Theorem~\ref{thm:main} can be checked by applying, e.g., the breadth-first search (BFS) algorithm. Specifically, one can augment the graph $G$ by adding a new node $\gamma$ and $m$ new edges $(\gamma,\beta_1),\ldots, (\gamma,\beta_m)$. Using BFS with $\gamma$ the starting vertex, we check whether every node of $V^\alpha$ is reachable from $\gamma$. Thus, the time complexity for checking item~1 is $O(|V| + |E|)$.   
For the remainder of the subsection, we exhibit a polynomial algorithm for checking the second item of Theorem~\ref{thm:main}. We do so by relating the condition~\eqref{eq:core} to the solution of a relevant maximum flow problem. 

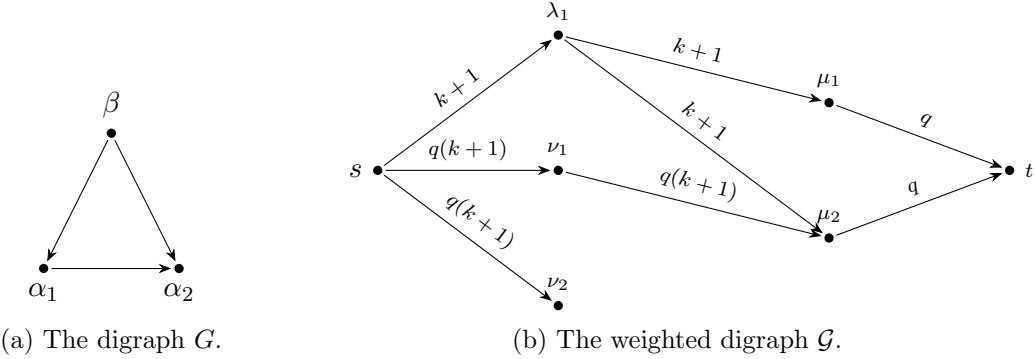
\begin{figure*}[ht]
  \centering
  %––––– Left subfigure –––––%
  \begin{subfigure}[t]{0.25\textwidth}
    \centering
    \begin{tikzpicture}[
    >=Stealth,            % Stealth arrowheads
    shorten >=1pt,        % gap between arrowheads and nodes
    scale=1.8,
    % every node/.style={font=1 pt},
    shorten <=1pt
  ]

  %––––––––––––––––––––––––––––––––––––––––––––––––––––%
  % Original G on the left
  %––––––––––––––––––––––––––––––––––––––––––––––––––––%
  \coordinate (P) at (-1.2,-0.6);
  % \node[above] at (P) {\(G:\)};

  \node[
    circle,fill=black,inner sep=0pt,minimum size=3.5pt,
    label=above:{\(\beta\)}
  ] (b)  at (0,0)   {};
  \node[
    circle,fill=black,inner sep=0pt,minimum size=3.5pt,
    label=below:{\(\alpha_1\)}
  ] (a1) at (-0.5,-1) {};
  \node[
    circle,fill=black,inner sep=0pt,minimum size=3.5pt,
    label=below:{\(\alpha_2\)}
  ] (a2) at (0.5,-1) {};

  \draw[->] (b)  -- (a1);
  \draw[->] (b)  -- (a2);
  \draw[->] (a1) -- (a2);
  
\end{tikzpicture}
    \caption{The digraph $G$.}
    \label{fig:small_graph}
  \end{subfigure}
  \hspace{0.02\textwidth}  % add horizontal spacing to center as a pair
  %––––– Right subfigure –––––%
  \begin{subfigure}[t]{0.55\textwidth}
    \centering
    \begin{tikzpicture}[
    >=Stealth,            % Stealth arrowheads
    shorten >=1pt,
    shorten <=1pt,
    scale = 1.2,
    every node/.style={color=black},  % all labels black
    every path/.style={draw=black}    % all lines black
  ]

  %––––––––––––––––––––––––––––––––––––––––––––––––––––%
  % adjustable spacing parameters
  %––––––––––––––––––––––––––––––––––––––––––––––––––––%
  \def\dxone{2cm}       % between col1 and col2
  \def\dxtwo{3cm}       % between col2 and col3
  \def\dxthree{2cm}     % between col3 and col4
  \def\dy{1.5cm}        % vertical spacing unit

  %% Column 1: single node s
  \node[circle, fill, inner sep=0pt, minimum size=3.5pt, label=left:{$s$}]
    (s) at (0,0) {};

  %% Column 2: 3 nodes (λ₁, ν₁, ν₂)
  \foreach [count=\i] \name in {\lambda_1,\nu_1,\nu_2}
  {
    \pgfmathsetmacro{\y}{(2-\i)*1} 
    \node[circle, fill, inner sep=0pt, minimum size=3.5pt,
          label=above:{\scriptsize \(\name\)}]
      (c2\i) at (\dxone,{(\y)*\dy}) {};
  }

  %% Column 3: 2 nodes (μ₁, μ₂)
  \foreach [count=\j] \rname in {\mu_1,\mu_2}
  {
    \pgfmathsetmacro{\y}{(1.5-\j)*1}
    \node[circle, fill, inner sep=0pt, minimum size=3.5pt,
          label=above:{\scriptsize \(\rname\)}]
      (c3\j) at ({\dxone+\dxtwo},{(\y)*\dy}) {};
  }

  %% Column 4: single node t
  \node[circle, fill, inner sep=0pt, minimum size=3.5pt, label=right:{\scriptsize $t$}]
    (t) at ({\dxone+\dxtwo+\dxthree},0) {};

  %% Edges (all black, uniform \scriptsize labels, sloped)
  \draw[->] (s)   -- node[midway, above, sloped]{\scriptsize ${k}$} (c21);
  \draw[->] (s)   -- node[midway, above, sloped]{\scriptsize $qk$} (c22);
  \draw[->] (s)   -- node[midway, above, sloped]{\scriptsize $qk$} (c23);

  \draw[->] (c21) -- node[midway, above, sloped]{\scriptsize ${k}$} (c31);
  \draw[->] (c21) -- node[midway, above, sloped]{\scriptsize ${k}$} (c32);

  \draw[->] (c22) -- node[midway, above, sloped]{\scriptsize $qk$} (c32);

  \draw[->] (c31) -- node[midway, above, sloped]{\scriptsize $q$} (t);
  \draw[->] (c32) -- node[midway, above, sloped]{\scriptsize $q$} (t);

\end{tikzpicture}
    \caption{The weighted digraph $\mathcal{G}$.}
    \label{fig:small_flow}
  \end{subfigure}

  %––––– Global caption –––––%
  \caption{Illustration of the graph $G$ (left) and the corresponding flow graph $\mathcal{G}$ (right).}
  \label{fig:combined_graphs}
\end{figure*}

Given $G = (V, E)$ and $k,q\in \N$, we construct a weighted digraph $\mathcal{G}=(\mathcal{V}, \mathcal{E}, c)$, with $c: \mathcal{E}\to \R_{\geq 0}$ specifying the edge weights, as follows: 
\begin{enumerate}
    \item The node set $\mathcal{V}$ is given by 
    $$
    \mathcal{V}=\{s\} \cup \mathcal{V}_L \cup \mathcal{V}_R \cup\{t\},
    $$
    where $\cV_L := \left\{\lambda_1, \ldots, \lambda_m, \nu_1, \ldots, \nu_n\right\}$ and $\cV_R := \left\{\mu_1, \ldots, \mu_n\right\}$. 

    \item The edge set $\cE$ is given by
    $$
    \cE = \cE_{sL} \cup \cE_{LR} \cup \cE_{Rt},
    $$
    where $\cE_{sL}$, $\cE_{LR}$, and $\cE_{Rt}$ are defined as follows:
    % \begin{equation*}
    % \left\{
    % \begin{aligned}
    % \mathcal{E}_{sL}&=\{(s,v)\mid v\in \cV_L\},
    % \\
    % \mathcal{E}_{LR}&=\left\{\left(\lambda_i, \mu_j\right) \mid\left(\lambda_i, \mu_j\right) \in E\right\} \cup\left\{\left(\nu_i, \mu_{j}\right) \mid\left(\mu_i, \mu_{j}\right) \in E\right\}, \\
    % \mathcal{E}_{Rt}&=\left\{\left(v, t\right) \mid v \in\cV_R\right\}.
    % \end{aligned}
    % \right.
    % \end{equation*}

    % \begin{equation*}
\begin{equation*}
\begin{aligned}
\mathcal{E}_{sL} &= \{(s,v)\mid v\in \cV_L\},\\
\mathcal{E}_{LR} &= 
    \{(\lambda_i, \mu_j)\mid (\lambda_i,\mu_j)\in E\}\cup\,
    \{(\nu_i, \mu_j)\mid (\mu_i,\mu_j)\in E\},\\
\mathcal{E}_{Rt} &= \{(v,t)\mid v\in \cV_R\}.
\end{aligned}
\end{equation*}

\item The weight $c:\cE\to \R_{\geq 0}$ is given by
\begin{subequations}
\label{eq:def_capacity_f}
\begin{align}
\mbox{On } \cE_{sL}:   c(s,\lambda_i):= & {k},   \mbox{ for } i \in [m], \\
 c(s, \nu_j) := & qk,  \mbox{ for } j \in [n], \\ 
\mbox{On } \cE_{LR}:    c(\lambda_i,\mu_j):= &{k},   \mbox{ for }
(\lambda_i,\mu_j) \in \cE_{LR}, \\
 c(\nu_i, \mu_j):= & qk,  \mbox{ for } (\nu_i,\mu_j) \in \cE_{LR}, \\
\mbox{On } \cE_{Rt}:   c(\mu_j, t):= & q,  \mbox{ for } j \in [n].
\end{align}
\end{subequations}
\end{enumerate}
See Figure~\ref{fig:small_flow} for an illustration.

% {  use subfigure for the $G$}

It is clear that $\mathcal{G}$ has $(2n + m + 2)$ nodes and $(2n + m + |E|)$ edges, so the time complexity for constructing $\mathcal{G}$ is $O(|V| + |E|)$. 

We now state the third main result of the paper:   

\begin{Theorem}
\label{thm:main2}
Given $G$ and $(k, q)$, let the weighted digraph $\mathcal{G}$ be given as above. 
Let $\cF$ be the set of flows on $\cG$, where $s$ is the source,  $t$ is the sink, and the edge weights are the capacities. 
Consider the maximum flow problem: 
\begin{equation}
\theta(k,q) := \max_{f\in \cF} |f|.
\label{prob:flow1} \tag{$P_1$}
\end{equation}
Then, condition~\eqref{eq:core} holds if and only if the solution to the above problem is $nq$, i.e., $\theta(k,q)=nq$.
\end{Theorem}

The class of maximum flow problems can be solved efficiently in polynomial time. See, for example, algorithms~\cite{karzanov1974determining,goldberg1985new,tarjan1984simple} with complexity $O(|V|^3)$ and algorithm~\cite{cherkassky1977algorithm} with complexity $O(\sqrt{|E|}|V|^2)$, all of which are {\it independent} of $k$ and $q$. 

It is also clear that the time complexity of the maximum flow problem~\eqref{prob:flow1} dominates over all the other steps (i.e., checking item~1 and construction of the weighted digraph $\mathcal{G}$) for checking the necessary and sufficient condition in Theorem~\ref{thm:main}. We summarize the above discussions in Algorithm~\ref{alg:maxflow}.

%We provide the proof in section~\ref{sec:proof_main2}.

Finally, we note that for the special case $q = 1$, the authors outlined in~\cite{liu2013structural} an algorithm by relating their necessary and sufficient condition to the solution of a maximal matching problem defined on a larger graph whose order is $O(k|V|)$ and whose size is $O(k|E|)$. Using the algorithm in~\cite{micali1980maximum}, the time complexity of their algorithm is then $O(\sqrt{k|V|} k|E|)$. 

\begin{algorithm}[h]
\DontPrintSemicolon      % optional: suppress automatic semicolons
% \textbf{Initialization:} sparsity pattern $\mathbb{S}$
\KwIn{sparsity pattern $\mathbb{S}$, $(k,q)\in \N\times \N$}
use BFS to check item~1 of Theorem~\ref{thm:main} \tcp*[r]{BFS: $O(|V|+|E|)$}
\If{item~1 of Theorem~\ref{thm:main} is not satisfied}
{\textbf{return }\textbf{False}}
\Else{
{construct the weighted digraph $\mathcal{G}$ \tcp*[r]{graph construction: $O(|V|+|E|)$}
solve Problem~\eqref{prob:flow1} for $(k,q)$ \tcp*[r]{maximum flow: $O(|V|^3)$}
\If{$\theta(k,q) = nq$ 
}
{\textbf{return }\textbf{True}} 
\Else{
\textbf{return }\textbf{False}}
}
}
\caption{Check whether $G$ is structurally controllable for a given $(k,q)$}
\label{alg:maxflow}
\end{algorithm}

\subsubsection{Algorithm for Proposition~\ref{prop:KS}}

As mentioned above, if $\mathcal{K}(\mathbb{S})$ is nonempty, then $k^*$ is bounded above by $n$. In particular, $\theta(n, q)=nq$, for all $q\in \N$. 

{  
Conversely, if $\mathcal{K}(\mathbb{S})$ is empty, then there must exist a $q\in \N$ such that~\eqref{eq:core} is not satisfied. 
For this to hold, it is necessary that there exists a $V'\subseteq V^\alpha$ such that $\iN^\alpha(V') = \varnothing$. This, in particular, implies that if 
\begin{equation}\label{eq:minimumq}
q\geq \frac{k |\iN^\beta(V')|}{|V'|}+1,  
\end{equation}
then $k |\iN^\beta(V')| < q |V'|$ (so~\eqref{eq:core} is not satisfied).  
Since $|\iN^\beta(V')| \leq m$, the right hand side of~\eqref{eq:minimumq} is bounded above by $(mk + 1)$. In a nutshell, if $\mathcal{K}(\mathbb{S}) = \varnothing$, then the minimum value of $q$ such that~\eqref{eq:core} does not hold with $k \leq n$ is bounded above by $(mn + 1)$.    

% In fact, it is not hard to see that any such $q$ is bounded above $(mn + 1)$. 
% let $V' \subseteq V^\alpha$ be such that
% $$
% n |\iN^\beta(V')| + nq |\iN^\alpha(V')| < q |V'|,
% $$
% which implies that $n|\iN^\alpha(V')| < |V'|$. It follows that
% $$
% q \geq  \frac{n |\iN^\beta(V')|}{|V'| - n|\iN^\alpha(V')| } + 1 \geq mn + 1,
% $$
% where we use $|\iN^\beta(V')| \leq m$ to establish the second inequality. 

In fact, the value $q = mn + 1$ plays an essential role in determining the set $\mathcal{K}(\mathbb{S})$ and hence, computing the smallest integer $k^*$ in that set (if it is nonempty).  
Specifically, we have the following result:} 

% We show below that one can always choose $q = mn + 1$ (regardless of the sparsity pattern).  
% Moreover, we show that the integer $k^*$ coincides with the smallest integer $k$ such that { $\theta(k,mn + 1)=n(mn+1)$}. 

% Precisely, we establish the following result: 

\begin{proposition}
Suppose that the sparsity pattern $\mathbb{S}$ satisfies item~1 of Theorem~\ref{thm:main}; { then, an integer $k\in \N$, with $1\leq k \leq n$, belongs to $\mathcal{K}(\mathbb{S})$ if and only if $\theta(k, mn+1) = n(mn+1)$.}
\end{proposition}
% { {proposition statement revised}}

\begin{proof}
First, note that if $k\in \mathcal{K}(\mathbb{S})$, then by~\eqref{eq:KS} and by item~2 of Theorem~\ref{thm:main}, we have that~\eqref{eq:core} is satisfied for the given $k$ and for all $q\in \N$, which include the case $q = mn + 1$. { By Theorem~\ref{thm:main2}, $\theta(k,mn+1)=n(mn+1)$.

We now show that if $\theta(k,mn+1)=n(mn+1)$, then $k\in \mathcal{K}(\mathbb{S})$.  }
%Assume that~\eqref{prob:flow1} has a solution of $nq$ for some $0\leq k\leq n-1$ and $q=mn+1$.
By Theorem~\ref{thm:main2}, we have that~\eqref{eq:core} is satisfied for $(k,q) = (k, mn+1)$:
\begin{equation}
\label{eq:cond_q1}
k |\iN^\beta(V')| + k(mn+1) |\iN^\alpha(V')|  \geq \\ (mn+1)|V'|, \quad \mbox{for all }V'\subseteq V^\alpha.
\end{equation}
We claim that %$\iN^\alpha(V')\neq \varnothing$ in the above inequality, and
\begin{equation}\label{eq:diffleq0}
|V'|-k|\iN^\alpha(V')|\leq 0, \quad \mbox{for all } V'\subseteq V^\alpha.
\end{equation} 
To wit, we assume to the contrary that there exists a $V'\subseteq V^\alpha$ such that $|V'|-k|\iN^\alpha(V')|> 0$. But then, using~\eqref{eq:cond_q1} and the fact that $k\leq n$, we obtain that 
$$mn + 1\leq \frac{k|\iN^\beta(V')|}{|V'|-k|\iN^\alpha(V')|}\leq km<mn+1,$$
which is a contradiction. We have thus established the claim. It then follows from~\eqref{eq:diffleq0} that for all $q\in \N$, 
\begin{equation}
\label{eq:cond_q2}
k |\iN^\beta(V')|  \geq 0 \geq q (|V'|-k|\iN^\alpha(V')|),  \mbox{ for all }V'\subseteq V^\alpha,
\end{equation}
so~\eqref{eq:core} holds for the given $k$ and for all $q\in \N$, which implies that $k\in \mathcal{K}(\mathbb{S})$. 
\end{proof}

\begin{algorithm}[h]
\DontPrintSemicolon      % optional: suppress automatic semicolons
% \textbf{Initialization:} sparsity pattern $\mathbb{S}$
\KwIn{sparsity pattern $\mathbb{S}$}
use BFS to check item~1 of Theorem~\ref{thm:main} \tcp*[r]{BFS: $O(|V|+|E|)$}
\If{item~1 of Theorem~\ref{thm:main} is not satisfied}
{\textbf{return }$k^*= \infty$ 
} 
solve Problem~\eqref{prob:flow1} for $(n,mn+1)$ \tcp*[r]{maximum flow: $O(|V|^3)$}
\If{$\theta(n,mn+1) \neq { n(mn+1)}$}
{\textbf{return }$k^*= \infty$ } 
\Else{
    $k_{\min} \gets 1$, \quad $k_{\max} \gets n$ 
    
    \While{$k_{\min} < k_{\max}$}{
        $k \gets \left\lfloor \frac{k_{\min} + k_{\max}}{2} \right\rfloor$ 
        solve Problem~\eqref{prob:flow1} for $(k,mn+1)$ 
        
        \If{$\theta(k,mn+1)={ n(mn+1)}$}{
            $k_{\max} \gets k$ 
        }
        \Else{
            $k_{\min} \gets k + 1$ 
        }
    }
    $k^* \gets k_{\min}$ \tcp*[r]{Binary Search: $O(|V|^3\log(|V|))$}
    {\textbf{return }$k^*= k_{\min}$}
}
\KwOut{$k^*$}
\caption{Compute $k^*$}
\label{alg:k^*}
\end{algorithm}

With the result above, we now present an algorithm for computing the value of $k^*$.  
The algorithm comprises two parts: The first part checks whether $\mathcal{K}(\mathbb{S})$ is nonempty. If $\mathcal{K}(\mathbb{S})$ is nonempty, then the second part computes the exact value of~$k^*$. 

\xc{Part 1: Check nonemptiness of $\mathcal{K}(\mathbb{S})$.} We use  BFS to check whether item~1 of Theorem~\ref{thm:main} is satisfied and use the maximum flow algorithm to check whether the solution { $\theta(n, mn+1)$ to Problem~\eqref{prob:flow1} is { $n(mn+1)$}. If either one of the two conditions is not satisfied, then $k^*= \infty$ and $\mathcal{K}(\mathbb{S}) = \varnothing$. }

\xc{Part 2: Compute the value of $k^*$.}
If $k^* < \infty$, we then proceed
to computing its exact value. Using binary search, we find the smallest $k$ such that the solution $\theta(k, mn + 1)$ to Problem~\eqref{prob:flow1} is { $n(mn+1)$}. 

The complete procedure is summarized in Algorithm~\ref{alg:k^*}.

\section{Proof of Theorem~\ref{thm:main}}
\label{sec:proof_main}

Note that if $G$ is structurally controllable for some $(k,q)$, then it is necessary that for every node $\alpha_i$, there exists at least a node $\beta_j$ and a path from $\beta_j$ to $\alpha_i$. This holds because otherwise, one can relabel the $\alpha$-nodes (if necessary) in a way such that any matrix pair $(A', B')$, with $\left[A' \  B'\right] \in \bS$, is in the form:
\[
A' = 
\begin{bmatrix}
A'_{11} & A'_{12} \\
0       & A'_{22}
\end{bmatrix},
\quad
B' =
\begin{bmatrix}
B'_1 \\
0
\end{bmatrix}.
\]
Now, consider the ensemble system~\eqref{eq:feltv}, with 
\(
\left[A_p[\ell] \ B_p[\ell]\right] \in \bS
\)
for all \( p = 1,\ldots,q \) and \( \ell = 1,\ldots,k \).
We decompose the state of the individual system \( x_p(t) = \bigl(x_{p,1}(t), x_{p,2}(t)\bigr) \) in accordance with the decomposition of the pair \(\bigl(A_p, B_p\bigr)\). Then, the dynamics of \( x_p(t) \) can be expressed as
$$
\begin{bmatrix}
\dot{x}_{p,1}(t) \\
\dot{x}_{p,2}(t)
\end{bmatrix}
= 
\begin{bmatrix}
A_{p,11}[\sigma(t)] & A_{p,12}[\sigma(t)] \\
0                  & A_{p,22}[\sigma(t)]
\end{bmatrix}
\begin{bmatrix}
x_{p,1}(t) \\
x_{p,2}(t)
\end{bmatrix}
\\
+ 
\begin{bmatrix}
B_{p,1}[\sigma(t)] \\
0
\end{bmatrix}
u(t). 
$$
It should be clear that the sub-state \( x_{p,2}(t) \) is not controllable. Thus, any individual system and hence, the ensemble system is not controllable.

The above arguments imply that the first item of Theorem~1 is necessary for \(G\) to be structurally controllable for \((k,q)\), for any \((k,q) \in \mathbb{N} \times \mathbb{N}\). In the sequel, we assume that this item is satisfied and show that the second item is necessary and sufficient for \(G\) to be structurally controllable for a given pair \((k,q)\).

Given the sparsity pattern $\bS$ and the integer $q$, we let 
$\hat{\bS} \subset \mathbb{R}^{nq \times (nq + m)}$ 
be the sparsity pattern associated with the $q$-ensemble system. 
Precisely, we define:
\begin{equation}
\label{eq:bS}
\hat{\bS}
:= 
\left\{
\begin{bmatrix}
A_{1}& &  & B_{1} \\
&  \ddots & &\vdots \\
&  &A_{q} & B_{q}
\end{bmatrix}
 \middle| 
[A_{p} \ B_{p}] \in \bS 
%\text{ for all } 
, \ \forall p = 1,\ldots,q
\right\}.
\end{equation}

%  \[
%     \hat{\bS} =
%     \left[
%     \begin{array}{ccccccccccccccc:cc}
%     0 & \star & 0 & 0 & 0 & & & & & & & & & &      &\star &\star\\
%     0 & 0     & 0 & 0 & 0 & & & & & & & & & &      & 0    &\star\\             \star & 0 & 0 & 0 & 0 & & & & & & & & & &      & 0    & 0   \\
%     0 & \star & \star & 0 & 0 & & & & & & & & & &  & 0    & 0   \\
%     0 & 0 & \star & \star & 0 & & & & & & & & & &  & 0 & 0      \\
%     \hline
%     & & & & & 0 & \star & 0 & 0 & 0 & & & & &      &\star &\star \\
%     & & & & & 0 & 0     & 0 & 0 & 0 & & & & &      &0 &\star     \\
%     & & & & & \star & 0 & 0 & 0 & 0 & & & & &      &0 &0          \\
%     & & & & & 0 & \star & \star & 0 & 0 & & & & &  &0 &0          \\
%     & & & & & 0 & 0 & \star & \star & 0 & & & & &  &0 &0           \\
%     \hline
%     & & & & & & & & & &  0 & \star & 0 & 0 & 0      &\star &\star     \\
%     & & & & & & & & & &  0 & 0     & 0 & 0 & 0      &0     &\star     \\
%     & & & & & & & & & &  \star & 0 & 0 & 0 & 0      &0     &0         \\
%     & & & & & & & & & &  0 & \star & \star & 0 & 0  &0     &0         \\
%     & & & & & & & & & &  0 & 0 & \star & \star & 0  &0     &0          \\
%     \end{array}
%     \right]
% \]

Let $\hat{G} = (\hat{V}, \hat{E})$ be the digraph associated with $\hat{\bS}$, 
and let $\hat V^\alpha$ and $\hat V^\beta$ be similarly defined.
Naturally, we can decompose 
\[
\hat{V}^\alpha  =  \cup_{p=1}^{q}\hat V^\alpha_p,
\]
such that the subgraph $\hat{G}_p$ of $\hat{G}$ induced by 
$\hat V^\beta \,\cup\, \hat V^\alpha_p$ corresponds to the $p$th individual system 
(note that all the $G_p$'s are isomorphic to $G$). More specifically, we write
$$
\hat{V}^\beta  = \{\,\hat{\beta}_1,\dots,\hat{\beta}_m\}
\quad\text{and}\quad
\hat{V}^\alpha_p  = \{\,\hat{\alpha}_{p,1},\dots,\hat{\alpha}_{p,n}\},
\quad \text{for } p = 1,\dots,q.
$$

% { 
The nodes are labeled such that the edge set $\hat{E}_p$ of $\hat{G}_p$ can be expressed as
$$
\hat{E}_p  :=  
\{\,(\hat{\alpha}_{p,i},\hat{\alpha}_{p,j}) \mid (\alpha_i ,\alpha_j) \in E\}  \cup   
\{\,(\hat{\beta}_i ,\hat{\alpha}_{p,j}) \mid (\beta_i,\alpha_j) \in E\}.
$$
% }

It should be clear that the edge set $\hat{E}$ of $\hat{G}$ is the union of the $\hat{E}_p$'s.
Denoted by $\pi_p : G \to \hat{G}_p$ the isomorphism, which sends
$\alpha_i$ (resp.\ $\beta_i$) to $\hat{\alpha}_{p,i}$ (resp.\ $\hat{\beta}_i$).
We need the following result, which is an immediate consequence of the above construction:

\begin{Lemma}\label{lemma:2}

For any $V' \subseteq V^\alpha$ and for any $p', p \in \{1,\dots,q\}$,
\(
N^\beta_{\mathrm{in}}\bigl(\pi_p(V')\bigr)
 = 
N^\beta_{\mathrm{in}}\bigl(\pi_{p'}(V')\bigr).
\)
\end{Lemma}

% {  
% It is clear that $G$ is structurally controllable for $(k,q)$ if and only if $\hat G$ is for $(k,1)$. 
% By Definition~\ref{def:structrl}, structural controllability for $(k,q)$ is defined through controllability of the stacked system~\eqref{eq:kswitchsys}. Since $\hat \bS $ in~\eqref{eq:bS} is the sparsity pattern of this stacked system, $G$ is structurally controllable for $(k,q)$ if and only if $\hat G$ is structurally controllable for $(k,1)$.
% }
It should be clear that $G$ is structurally controllable for $(k,q)$ if and only if $\hat{G}$ is for $(k,1)$. 
{ The following result can be adapted from~\cite{liu2013structural}, which provides a necessary and sufficient condition for $\hat G$ to be structurally controllable for $(k,1)$.}

\begin{Lemma}\label{lemma:3}
The sparsity pattern $\hat{\mathbb{S}}$ is structurally controllable for $(k,1)$ 
if and only if the following conditions hold:

\begin{enumerate}
    \item For every node $\hat{\alpha}_{p,i}$, there exist a node $\hat{\beta}_j$ and a path from $\hat{\beta_j}$ to $\hat{\alpha}_{p,i}$.
    \item We have that
    % \vspace{-3pt}
   \begin{equation}
    \label{eq:12}
k\,\lvert N_{\mathrm{in}}(\hat V')\rvert  \ge  \lvert \hat V'\rvert,
\quad
\text{for any } \hat V' \subseteq \hat V^\alpha.
\end{equation}
\end{enumerate}
\end{Lemma}

{ 
\begin{Remark}\normalfont
We omit the proof of the above lemma as it is essentially the same as
Theorem~9 of~\cite{liu2013structural}, specialized to the case where the sparsity patterns of the subsystems are the same, given by $\hat \bS$, with $\hat G$ the graph representation. Consequently, item~2 of Theorem~9 of~\cite{liu2013structural}, with their definition of $S$-dilation given in Definition~12 of the paper, is reduced to our condition~\eqref{eq:12}. \hfill{\qed}
% applied to the switched system with $k$ subsystem digraphs, each identical to $\hat{\mathcal{G}}$. Here, for any $\hat{V}'\subseteq \hat{V}^\alpha, T_{\ell}(\hat{V}')$ denotes the set of vertices in the $\ell$-th subsystem digraph that have an edge into $\hat V'$, and $|T(\hat{V}')|=\sum_{\ell=1}^k |T_{\ell}(\hat{V}')|$. Since all subsystems share the same sparsity pattern, for every $\hat{V}' \subseteq \hat{V}^\alpha$ one has $T_{\ell}(\hat{V}')=N_{\text {in }}(\hat{V}')$ for all $\ell=1, \ldots, k$, and hence the $S$-dilation condition $\sum_{\ell=1}^k|T_{\ell}(\hat{V}')| \geq|\hat{V}'|$ becomes $k|N_{\text {in }}(\hat{V}')|\geq|\hat{V}'|$. 
% \hfill{$\qedsymbol$}
\end{Remark}
}

Now, to establish Theorem~\ref{thm:main}, it remains to show the following result:

\begin{proposition}
\label{prop:double_star}
    Condition~\eqref{eq:core} holds if and only if ~\eqref{eq:12} holds.
\end{proposition}

\begin{proof}
We first show that~\eqref{eq:12} implies~\eqref{eq:core} and then establish the converse.

\medskip
\noindent
\textit{Proof that~\eqref{eq:12} implies~\eqref{eq:core}.}
Given an arbitrary \(V' \subseteq V^\alpha\), let
\[
\hat V'_p  :=  \pi_p(V')
\quad\mbox{and}\quad
\hat V'  :=  \cup_{p=1}^{q} \hat V'_p.
\]
Since \(\hat{G}_p\) is isomorphic to \(G\), we have
\[
N_{\mathrm{in}}^\alpha(\hat V'_p)
 =  \pi_p\bigl(N_{\mathrm{in}}^\alpha(V')\bigr).
\]
It follows that
\[
\bigl\lvert N_{\mathrm{in}}^\alpha(\hat V'_p) \bigr\rvert
 = 
\bigl\lvert N_{\mathrm{in}}^\alpha(V') \bigr\rvert,
\quad \text{for all } p = 1,\dots,q,
\]
and hence,
\begin{equation}\label{eq:13}
\lvert N_{\mathrm{in}}^\alpha(\hat V')\rvert
 = 
q\,\bigl\lvert N_{\mathrm{in}}^\alpha(V')\bigr\rvert.
\end{equation}

Next, by Lemma~\ref{lemma:2}, we have
\[
N_{\mathrm{in}}^\beta(\hat V')
 = 
N_{\mathrm{in}}^\beta(\hat V'_p)
 = 
\pi_p\bigl(N_{\mathrm{in}}^\beta(V')\bigr),
\quad \text{for all } p = 1,\dots,q,
\]
and hence,
\begin{equation}\label{eq:14}
\bigl\lvert N_{\mathrm{in}}^\beta(\hat V')\bigr\rvert
 = 
\bigl\lvert N_{\mathrm{in}}^\beta(V')\bigr\rvert.
\end{equation}

Using~\eqref{eq:13} and~\eqref{eq:14}, we obtain
$$
k\lvert N_{\mathrm{in}}(\hat V')\rvert
=
k\Bigl[
\bigl\lvert N_{\mathrm{in}}^\beta(\hat V')\bigr\rvert
 + 
\bigl\lvert N_{\mathrm{in}}^\alpha(\hat V')\bigr\rvert
\Bigr]
=
k\bigl\lvert N_{\mathrm{in}}^\beta(V')\bigr\rvert
 + 
k\,q\,\bigl\lvert N_{\mathrm{in}}^\alpha(V')\bigr\rvert.
$$
Thus, assuming that~\eqref{eq:12} holds, we get
$$
k\bigl\lvert N_{\mathrm{in}}^\beta(V')\bigr\rvert
 + 
k\,q\,\bigl\lvert N_{\mathrm{in}}^\alpha(V')\bigr\rvert
 = 
k\lvert N_{\mathrm{in}}(\hat V')\rvert
  \geq 
\lvert \hat V'\rvert
 = 
\sum_{p=1}^q \lvert \hat V'_p\rvert
 = 
q\,\lvert V'\rvert.
$$
This completes the proof that~\eqref{eq:12} implies~\eqref{eq:core}. 

\medskip
\noindent
\textit{Proof that~\eqref{eq:core} implies~\eqref{eq:12}.} Now, let $\hat V' = \cup_{p = 1}^q \hat V'_p$ be an arbitrary subset of $\hat V^\alpha$, with $\hat V'_p \subseteq \hat V^\alpha_p$. We define 
$$V':= \cup_{p = 1}^q V'_p, \quad \mbox{with } V'_p := \pi_p^{-1}(\hat V'_p).$$
Note that $V'$ is a subset of $V^\alpha$. We have that  
\begin{equation}\label{eq:part3}
|\iN^\beta(\hat V')| =  |\iN^\beta(\cup_{p = 1}^q \hat V'_p)| = |\iN^\beta(\cup_{p = 1}^q V'_p)| = |\iN^\beta(V')|. 
\end{equation}
Since each $V'_p$ is a subset of $V'$, we have that $\iN^\beta(V') \supseteq \iN^\beta(V'_p)$ and hence, 
\begin{equation}\label{eq:part4}
    |\iN^\beta(V')| \geq | \iN^\beta(V'_p)|\quad \mbox{for any } p = 1,\ldots, q. 
\end{equation}
Next, let 
$$
p^* \in \arg\min_{p} \{k|\iN^\alpha(\hat V'_p)| - |\hat V'_p| \}. 
$$
Then, 
\begin{equation}
\begin{aligned}\label{eq:part5}
k |\iN^\alpha(\hat V')| - |\hat V'| 
= \sum_{p = 1}^q \left [ k|\iN^\alpha(\hat V'_p)| - |\hat V'_p|\right ]  
\geq & \ q \left [ k|\iN^\alpha(\hat V'_{p^*})| - |\hat V'_{p^*}| \right ] \\=& \ q \left[ k |\iN^\alpha(V'_{p^*})| - |V'_{p^*}| \right].
\end{aligned}
\end{equation}
Combining the above arguments, we obtain that
$$
\begin{aligned}
k|\iN(\hat V')| - |\hat V'| = &k|\iN^\beta(\hat V')| + k|\iN^\alpha(\hat V')| - |\hat V'| \\
 \geq & k|\iN^\beta(V'_{p^*})| + q \left [ k|\iN^\alpha(V'_{p^*})| - |V'_{p^*}| \right ] \geq 0,
\end{aligned}
$$
where we use~\eqref{eq:part3},~\eqref{eq:part4}, and~\eqref{eq:part5} to establish the first inequality and the second inequality follows from the hypothesis that~\eqref{eq:core} holds. This completes the proof. 
\end{proof}

\section{Proof of Theorem~\ref{thm:main2}}
\label{sec:proof_flow_equiv}
% To establish Theorem~\ref{thm:main2}, we need the following as preliminary results.
% % Theorem~\ref{thm:flow_equiv} 
% We first construct a graph $\hat \cG$ which we will later relate to $\cG$ from~\eqref{prob:flow1}. An example of how $\hat \cG$, and $\cG$ is constructed from $G$ is provided in Figure~\ref{fig: small_graph},~\ref{fig:small_flow}, and~\ref{fig:big_flow}.

\begin{figure*}[ht]
  \centering
  % includegraphics options: width, height, scale, angle, etc.
  \includegraphics[width=0.78\textwidth]{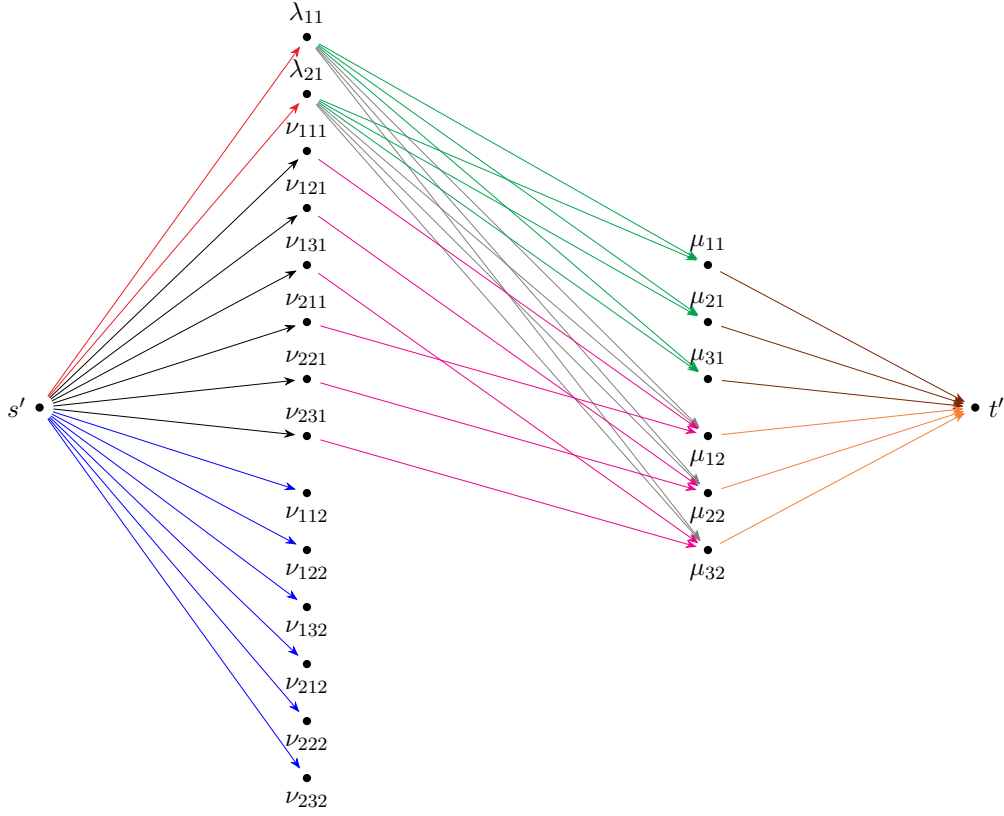}
  \caption{Illustration of the digraph $\hat{\mathcal{G}}$ for $G$ given in Figure~\ref{fig:small_graph} and $(k,q) = (2,3)$. 
Edges are colored according to their corresponding edges in the digraph $\mathcal{G}$ in
Fig.~\ref{fig:small_flow}. For any edge $(u,v)$ in $\mathcal{G}$, all edges in  $\phi^{-1}(u,v)$, with $\phi$ the graph homomorphism defined in~\eqref{eq:def_phi}, 
are plotted in the same color. Specifically, the red edges are those in $\phi^{-1}(s,\lambda_1)$, the blacks are in $\phi^{-1}(s,\nu_1)$, the blues are in $\phi^{-1}(s,\nu_2)$, the greens are in $\phi^{-1}(\lambda_1,\mu_1)$, the grays are in $\phi^{-1}(\lambda_1,\mu_2)$, the magentas are in $\phi^{-1}(\nu_1,\mu_2)$, the browns are in $\phi^{-1}(\mu_1,t)$, and the oranges are in $\phi^{-1}(\mu_2,t)$. For each color, the number of edges in $\hat{\mathcal{G}}$ is determined by the weight of the corresponding edge in $\mathcal{G}$. Every edge of $\hat{\mathcal{G}}$ has
weight~$1$.}

\label{fig:big_flow}
\end{figure*}

Given the digraph $G = (V, E)$ and the pair $(k,q)\in \N\times \N$, we introduce below another weighted digraph $\hat {\mathcal{G}}=(\hat {\mathcal{V}},\hat{\mathcal{E}},\hat{c})$: 
\begin{enumerate}
    \item The node set $\hat {\cV}$ is given by 
    $$
    \hat{\mathcal{V}}=\left\{s \right \} \cup \hat{\mathcal{V}}_L \cup \hat{\mathcal{V}}_R \cup\left \{t\right \},
    $$
    where $s$ and $t$ are the source and sink of $\hat G$, respectively, and $\hat \cV_L$ and $\hat \cV_R$ are defined as follows:
    \begin{align*}
\hat{\mathcal{V}}_L 
    &:= \{\lambda_{\ell i} 
            \mid \ell\in[k],\, i\in[m]\}
        \cup
        \{\nu_{\ell p j} 
            \mid \ell\in[k],\, p\in[q],\, j\in[n]\},
        \\
\hat{\mathcal{V}}_R
    &:= \{\mu_{p j} 
            \mid p\in[q],\, j\in[n]\}.
\end{align*}

    so $|\hat \cV_L| = k(m + nq)$ and $|\hat \cV_R| = nq$.

    \item The edge set $\hat{\cE}$ is given by
    $$
    \hat{\cE} = \hat{\cE}_{sL} \cup \hat{\cE}_{LR} \cup \hat{\cE}_{Rt},
    $$
    where $\hat{\cE}_{sL}$, $\hat{\cE}_{LR}$, and $\hat{\cE}_{Rt}$ are defined as follows: 
    \begin{equation}
    \label{eq:def_edge_set}
    \left \{
    \begin{aligned}
        \hat{\mathcal{E}}_{sL} := & \{(s, v)\mid v \in \hat \cV_L \},\\
        \hat{\mathcal{E}}_{LR}  := & \{(\lambda_{\ell i},\mu_{pj})\mid  \ell \in [k], p\in [q], \mbox{ and } (\lambda_i,\mu_j)\in E \} \, \\
        & \cup \{(\nu_{\ell pi},\mu_{pj})\mid \ell \in [k], p\in [q], 
        \mbox{ and } (\mu_i,\mu_{j})\in E \},\\
        \hat{\mathcal{E}}_{Rt} := & \{ \left (v, t \right) \mid v\in \hat \cV_R \}.
    \end{aligned}
    \right.
    \end{equation}
    \item The weight $\hat {c}:\hat{\cE}\to \R_{\geq 0}$ is defined as
    \begin{equation}
    \label{eq:def_capacity_f_bar}
    \hat{c}(e):=1, \quad \mbox{for all } e\in \hat{\mathcal{E}}.
    \end{equation}
\end{enumerate}
See Figure~\ref{fig:big_flow} for an illustration. 

We have the following result:

\begin{Theorem}
\label{thm:flow_equiv}
    Given $G$ and $(k,q)$, let $\hat{\mathcal{G}}$ be given as above.  Let $\hat \cF$ be the set of flows on $\hat{\mathcal{G}}$. Consider the following maximum flow problem~\ref{prob:flow2} on $\hat \cG$:
\begin{align}
    &\hat \theta(k,q) := \max_{\hat {f}\in \hat {\mathcal{F}}} |\hat {f}|.
    \label{prob:flow2} \tag{$P_2$}
\end{align}
Then, the solution to~\ref{prob:flow2} is the same as the solution to~\eqref{prob:flow1}, i.e., $\hat \theta(k,q)=\theta(k,q)$.  
\end{Theorem}

 % if the above holds, we prove the main theorem
It is well known~\cite{hall1998combinatorial} that the solution $\hat \theta(k,q)$ to~\eqref{prob:flow2} is the same as the solution to the maximal matching problem for the bipartite graph $(\hat{\mathcal{V}}_L,\hat{\mathcal{V}}_R, \hat{\mathcal{E}}_{LR} )$, which in turn decides whether~\eqref{eq:core} is satisfied. We elaborate below on this fact:

\hspace{.2cm}
{\it Proof of Theorem~\ref{thm:main2}:} 
For the given $G$ and $(k,q)$, let $\hat G$ be given as in Section~\ref{sec:proof_main}. By Proposition~\ref{prop:double_star},~\eqref{eq:core} holds if and only if~\eqref{eq:12} holds. 
Also, by Theorem~\ref{thm:flow_equiv}, $\theta(k,q) = nq$ if and only if $\hat \theta(k,q) = nq$. Thus, to establish Theorem~\ref{thm:main2}, it suffices to show that~\eqref{eq:12} holds if and only if $\hat \theta(k,q) = nq$. 
But it follows from~\cite{liu2013structural} (see Definition 11, Theorem 11, and the arguments in Section 3.3 of their paper) that~\eqref{eq:12} holds if and only if the bipartite graph $(\hat{\mathcal{V}}_L, \hat{\mathcal{V}}_R, \hat{\mathcal{E}}_{LR})$ has a right perfect matching, which holds if and only if $\hat \theta(k,q) = nq$. \hfill{\small $\blacksquare$}

    % We first show that the condition in~\eqref{eq:12} holds if and only if the solution $\hat \theta(k,q)$ to~\ref{prob:flow2} is $nq$. 
    % By  in~\cite[Theorem 11]{liu2013structural},~\eqref{eq:12} can be verified by checking the existence of a right perfect matching in the bipartite graph $(\hat{\mathcal{V}}_L, \hat{\mathcal{V}}_R, \hat{\mathcal{E}}_{LR})$. It is well known that such a matching can be found by solving a maximum flow problem~\cite{kleinberg2006algorithm}, which corresponds exactly to the flow problem we formulate in~\eqref{prob:flow2}. A right perfect matching exists if and only if the solution $\hat\theta(k,q)$ to~\eqref{prob:flow2} equals the cardinality of the right bipartite node set $nq$.

    % Combining Proposition~\ref{prop:double_star} and Theorem~\ref{thm:flow_equiv}, we obtain Theorem~\ref{thm:main2}. First, condition~\eqref{eq:12} holds if and only if condition~\eqref{eq:core} holds. These conditions are in turn equivalent to the existence of a flow of value $nq$ for Problem~\ref{prob:flow2}, which is further equivalent to Problem~\ref{prob:flow1} admitting a solution of value $nq$.

The remainder of the section is devoted to the proof of Theorem~\ref{thm:flow_equiv}, which is divided into two parts. In Subsections~\ref{ssec:equiflowpart1} and~\ref{ssec:equiflowpart2}, we establish the two inequalities $\theta(k,q) \geq \hat \theta(k,q)$ and $\theta(k,q) \leq \hat \theta(k,q)$, respectively. 

% \begin{align}
%     &\max_{f\in \cF} |f|
%     \label{prob:flow1} \tag{$P_1$}
%     \\
%     \mbox{s.t. }
%     &0 \leq f(v_i,v_j)\leq c(v_i,v_j),\forall (v_i,v_j)\in \mathcal{E},
%     \label{eq:capacity_con}
%     \\
%     &\sum_{(v_i,v_j) \in \mathcal{E}} f(v_i,v_j)=\sum_{(v_j,v_i) \in \mathcal{E}} f(v_j,v_i), \forall v_i\in \mathcal{V}_L\cup\mathcal{V}_R.
%     \label{eq:flow_balance}
% \end{align}

% \begin{Theorem}
% \label{thm:equiv_flow}
%     Condition~\eqref{eq:core} holds if and only if the solution to Problem \ref{prob:flow1} is $n$.
% \end{Theorem}

%\paragraph{Sketch of Proof:}
% We start by constructing the graph $\bar{\mathcal{G}}$ and then define the maximum flow problems~\ref{prob:flow2} on $\bar{\mathcal{G}}$.
%We prove the Theorem~\ref{thm:main2} in two steps. Firstly, we show that checking the condition in~\eqref{eq:12} for $(k,q)$ is equivalent to checking the solution to problem~\ref{prob:flow2}. In the next step, we show that if a solution of $n$ in problem~\ref{prob:flow1} exists, then~\ref{prob:flow2} has a feasible solution that achieves a maximum flow of $nq$. With this result, we can apply Proposition~\ref{prop:double_star} and are able to show that checking the condition in~\eqref{eq:core} can be done by solving for the solution in~\ref{prob:flow1}.

% \section{Proof of Theorem~\ref{thm:flow_equiv}}

%\subsection{Construction of $\hat{\mathcal{G}}$}

\subsection{Proof that $\theta(k,q) \geq \hat \theta(k,q)$}\label{ssec:equiflowpart1}

%{ {introduce homomorphism}}

% In this subsection, we establish the following proposition: 

% \begin{proposition}
% \label{prop:F_hat_to_F}
%     For any flow $\hat f \in \hat{\mathcal{F}}$, there exists a flow $f\in \mathcal{F}$ such that $|f| =  |\hat f|$.
% \end{proposition}
We show that for any flow $\hat f \in \hat{\mathcal{F}}$, there exists a flow $f\in \mathcal{F}$ such that $|f| =  |\hat f|$. 
The proof is constructive. 

To proceed, we first introduce the map $\phi:\hat{\mathcal{V}}\to \mathcal{V}$ defined as follows:
\begin{equation}
\label{eq:def_phi}
\phi(s):=s,\quad 
\phi(\lambda_{\ell i}):=\lambda_i,\quad
\phi(\nu_{\ell p j}):=\nu_j,\quad
\phi(\mu_{p j}):=\mu_{j},\quad
\phi(t):=t.
\end{equation}
The following lemma follows directly from the construction of $\cG$ and $\hat \cG$, and its proof is omitted. 

\begin{Lemma}
{ \label{lem:graph_homo}
    The map $\phi$ is a  graph homomorphism, i.e., if $(u,v)$ is an edge of $\hat \cG$, then $(\phi(u),\phi(v))$ is an edge of $\cG$.} 
\end{Lemma}

By Lemma~\ref{lem:graph_homo}, 
we can extend the domain of $\phi$ to include the edge set $\hat \cE$ of $\hat \cG$ and write $\phi: (\hat \cV, \hat \cE) \to (\cV, \cE)$. %{\color{red}(to be removed)Also, note that for any $e\in \cE$, we have that $\sum_{\hat e\in \phi^{-1}(e)} \hat c(\hat e) = c(e)$. }
Now, given the flow $\hat f\in \cF$, we define $f:  \cE\to \R_{\geq 0}$ as follows: 
\begin{equation}
    \label{eq:construct_F}
    f(e) :=\sum_{\hat{e}\in \phi^{-1}(e)}\hat{f}(\hat{e}), \quad \mbox{for all } e\in \mathcal{E}.
\end{equation}
We have the following lemma:

\begin{Lemma}
\label{lem:f_flow}
The map $f$ defined in~\eqref{eq:construct_F} is a flow. 
\end{Lemma}

\begin{proof} We will first show that $f$ satisfies the conservation of flow and then show that $f$ satisfies the capacity constraint. For ease of presentation, we extend the domain of $f$ (resp., $\hat f$) to $\cV\times \cV$ (resp., $\hat \cV\times \hat \cV$) by setting $f(u,v) := 0$ if $(u,v)\notin \cE$ (resp., $\hat f(u,v) := 0$ if $(u,v)\notin \hat \cE$).

\xc{Proof that $f$ satisfies the conservation of flow.}
    We show below that~\eqref{eq:flow_balance} is satisfied for every node in $\mathcal{V}_L\cup\mathcal{V}_R$: 
\begin{enumerate}
\item For any $\lambda_i\in \cV_L$, the only incoming neighbor is the source~$s$. Also, note that 
\begin{align}
\label{eq:phi_s_beta}
\phi^{-1}(s,\lambda_i) &= \{(s,\lambda_{\ell i}) \mid \ell \in [k]\}, \\
\phi^{-1}(\lambda_i,\mu_j) &= \{(\lambda_{\ell i},\mu_{pj}) \mid \ell\in [k], p\in [q]\}.\label{eq:phi_beta_alpha}
\end{align}
It then follows that
\noindent
\begin{equation}
\label{eq:verify_f_s_beta}
\begin{aligned}
f(s, \lambda_i)
= \sum_{\ell=1}^{k} \hat{f}(s, \lambda_{\ell i})
= &\sum_{\ell=1}^{k} \left [\sum_{p=1}^q \sum_{j=1}^n \hat{f}(\lambda_{\ell i}, \mu_{p j})\right]
\\
=&\sum_{j=1}^n\left [\sum_{\ell=1}^{k} \sum_{p=1}^q  \hat{f}\left(\lambda_{\ell i}, \mu_{p j}\right)\right ] 
=  \sum_{j=1}^n f\left(\lambda_i, \mu_{j}\right),
% f\left(s, \lambda_i\right)= \sum_{\ell=1}^{k+1} \hat{f}\left(s, \lambda_{\ell i}\right)\\
\end{aligned}
\end{equation}
where the first equality follows from~\eqref{eq:construct_F} and~\eqref{eq:phi_s_beta}, the second equality follows from the fact that $\hat f$ satisfies the conservation of flow at the node $\lambda_{\ell i}$ for any $\ell\in [k]$, and the last equality follows from~\eqref{eq:construct_F} and~\eqref{eq:phi_beta_alpha}. 
%accordingly. Therefore, the conservation of flow is satisfied at node $\lambda_i$, for all $i\in [m]$. 
\item Similarly, for any $\nu_i\in \cV_L$, the only incoming neighbor is the source~$s$ and, moreover,
\begin{align}
\label{eq:phi_s_gamma}
\phi^{-1}(s,\nu_i) & = \{(s,\nu_{\ell p i}) \mid \ell\in [k], p\in [q]\}, \\
\label{eq:phi_gamma_alpha}
\phi^{-1}(\nu_i,\mu_j) & = \{(\nu_{\ell p i},\mu_{pj})\mid \ell \in [k], p\in [q]\}.
\end{align}
Thus, 
\begin{align*}
    f\left(s, \nu_i\right)
    = \sum_{\ell=1}^{k} \sum_{p=1}^q \hat{f}\left(s, \nu_{\ell p i}\right)
    &=\sum_{\ell=1}^{k} \left[\sum_{p=1}^q \sum_{j=1}^n \hat{f}\left(\nu_{\ell p i}, \mu_{p j}\right)\right ]
    \\
    &=\sum_{j=1}^n\left[ \sum_{\ell=1}^{k} \sum_{p=1}^q  \hat{f}\left(\nu_{\ell p i}, \mu_{p j}\right)\right] \\
    & =\sum_{j=1}^n f\left(\nu_i, \mu_{j}\right),
\end{align*}
where the first equality follows from~\eqref{eq:construct_F} and~\eqref{eq:phi_s_gamma}, the second equality holds due to the conservation of flow 
 at node $\nu_{\ell p i}$, for any $\ell\in [k]$ and $p\in [q]$, and the last equation follows from~\eqref{eq:construct_F} and~\eqref{eq:phi_gamma_alpha}. 

\item Finally, for any $\mu_j\in \cV_R$,
the only outgoing neighbor is the sink~$t$ and 
\begin{align}
\label{eq:phi_alpha_t}
\phi^{-1}(\mu_j, t) & = \{(\mu_{p j},t)\mid p\in [q]\}.
\end{align}
We then obtain that 
\begin{align*}
    f\left(\mu_{j}, t\right) & =\sum_{p=1}^q \hat{f}\left(\mu_{p j}, t\right)
    \\
    & =\sum_{p=1}^q\left[\sum_{i=1}^m  \sum_{\ell=1}^{k} \hat{f}\left(\lambda_{\ell i}, \mu_{p j}\right) \right.  \left. +\sum_{i=1}^n \sum_{\ell=1}^{k}  \hat{f}\left(\nu_{\ell p i}, \mu_{p j}\right)\right] \\
    & =\sum_{i=1}^m  \left [ \sum_{\ell=1}^{k} \sum_{p=1}^q \hat{f}\left(\lambda_{\ell i}, \mu_{p j}\right) \right ]+\sum_{i=1}^n \left [ \sum_{\ell=1}^{k} \sum_{p=1}^q \hat{f}\left(\nu_{\ell p i}, \mu_{p j}\right) \right ] \\
    & =\sum_{i=1}^m f\left(\lambda_i, \mu_{j}\right)+\sum_{i=1}^n f\left(\nu_i, \mu_{j}\right),
\end{align*}
where the first equality follows from~\eqref{eq:construct_F} and~\eqref{eq:phi_alpha_t}, the second equality follows from the  conservation of flow at node $\mu_{pj}$, for any $p\in [q]$, and the last equality follows from~\eqref{eq:construct_F},~\eqref{eq:phi_beta_alpha} and~\eqref{eq:phi_gamma_alpha}.
%Thus, we have completed the verification for every node in $\cV_L\cup\cV_R$ and this shows that $f$ satisfies the flow of conservation property.
\end{enumerate}

\xc{Proof that $f$ satisfies the capacity constraint.}
    %We proceed with the verification that for any given $\hat{F}$ in $\bar{\mathcal{G}}$, each element in the constructed flow $F$ does not exceed the capacity assigned in~\eqref{eq:def_capacity_f}. 
    %First of all, by the feasibility of $\hat{f}$, $\hat{f}$ should not exceed the capacity assigned in~\eqref{eq:def_capacity_f_bar}. The followings are true under~\eqref{eq:construct_F}:
We show below that $f$ satisfies $f(e) \leq c(e)$ for all $e\in \cE$, with $c$ defined in~\eqref{eq:def_capacity_f}. First, using~\eqref{eq:construct_F} and~\eqref{eq:phi_s_beta}, and the fact that $\hat f(\hat e) \leq 1$ for all $\hat e\in \hat \cE$, we obtain that
\begin{equation}
\label{eq:f_s_beta}
f\left(s, \lambda_i\right) = \sum_{\ell=1}^{k} \hat{f}\left(s, \lambda_{\ell i}\right)\leq {k}. 
\end{equation}
Using the same arguments, but with~\eqref{eq:phi_s_beta} replaced by~\eqref{eq:phi_s_gamma},~\eqref{eq:phi_gamma_alpha}, and~\eqref{eq:phi_alpha_t},  respectively, we have that
\begin{align*}
f\left(s, \nu_j\right)&= \sum_{\ell=1}^{k} \sum_{p=1}^q \hat{f}\left(s, \nu_{\ell p i}\right)\leq qk,\\
f\left(\nu_i, \mu_{j}\right)&= \sum_{\ell=1}^{k} \sum_{p=1}^q \hat{f}\left(\nu_{\ell p i}, \mu_{p j}\right)\leq qk,
\\
f\left(\mu_{j}, t\right)&= \sum_{p=1}^q \hat{f}\left(\mu_{p j}, t\right)\leq q.
\end{align*}
Finally, for the edges $(\lambda_i,\mu_j)\in \cE$, we have that
\begin{equation}
    f\left(\lambda_i, \mu_{j}\right)\leq f\left(s, \lambda_{i}\right) \leq {k},
\end{equation}
where the first inequality follows from~\eqref{eq:verify_f_s_beta} and the second inequality follows from~\eqref{eq:f_s_beta}.
%it is guaranteed that the flow going out of node $\lambda_i$ does exceed the one going into the node, which is $f(s,\lambda_i )$. This results in $$f\left(\lambda_i, \mu_{j}\right)\leq {k+1}$$ for~\eqref{eq:modify_capacity}. Consequently, after a direct comparison with the capacity assignment in~\eqref{eq:def_capacity_f}, we have verified the capacity constraint for $f$ in~\eqref{eq:construct_F}. Thus, $f$ is indeed a flow.
\end{proof}

With the lemmas above, we can now show that $\theta(k, q) \geq \hat{\theta}(k, q)$. 
Let $f$ be given as in~\eqref{eq:construct_F}. By Lemma~\ref{lem:f_flow}, $f\in \mathcal{F}$. Moreover,  
\begin{equation*}
    |f|=\sum_{j=1}^n f(\mu_j,t)=\sum_{p=1}^q\sum_{j=1}^n \hat f(\mu_{pj},t)=|\hat f|,
\end{equation*}
where the first and the last equality follow from the definition~\ref{eq:value}, and the second equality follows from~\eqref{eq:construct_F} and~\eqref{eq:phi_alpha_t}. 
\unskip\nobreak\hfill\mbox{\qed}

\subsection{Proof that $\theta(k,q)\leq \hat \theta(k,q)$}\label{ssec:equiflowpart2}

Similarly, in this subsection, we show that for any flow $f\in\mathcal{F}$, there exists a flow $\hat{f}\in \hat{\mathcal{F}}$ such that $|\hat f| =  |f|$. 
To this end, we consider the following map $\hat f: \hat \cE\to\R_{\geq 0}$ as follows:
\begin{subequations}\label{eq:construct_f_hat}
\begin{numcases}{\hat f(\hat e):=}
  \frac{1}{k}\,f\!\bigl(\phi(\hat e)\bigr)
    & if $\hat e\in\hat{\mathcal E}_{sL}$, \label{eq:f_hat_s_L}\\
  \frac{1}{kq}\,f\!\bigl(\phi(\hat e)\bigr)
    & if $\hat e\in\hat{\mathcal E}_{LR}$, \label{eq:f_hat_L_R}\\
  \frac{1}{q}\,f\!\bigl(\phi(\hat e)\bigr)
    & if $\hat e\in\hat{\mathcal E}_{Rt}$. \label{eq:f_hat_R_t}
\end{numcases}
\end{subequations}

% \begin{subequations}
% \label{eq:construct_f_hat}
% \begin{align}
% \label{eq:f_hat_s_L}
% \hat f (e) &= \frac{1}{k+1}f(\phi(e)),\mbox{ for all } e\in \hat \cE_{sL},\\
% \label{eq:f_hat_L_R}
% \hat f (e) &= \frac{1}{(k+1)q}f(\phi(e)),\mbox{ for all } e\in \hat \cE_{LR},\\
% \label{eq:f_hat_R_t}
% \hat f (e) &= \frac{1}{q}f(\phi(e)),\mbox{ for all } e\in \hat \cE_{Rt}.
% \end{align}
% \begin{align}
%     & \hat{f}\left(s, \lambda_{\ell i}\right)=\frac{1}{k+1} f\left(s, \lambda_i\right),\mbox{ for all }\ell\in [k+1]\mbox{ and }i\in [m],
%     \label{eq:construct_f'_s_beta}\\
%     & \hat{f}\left(s, \nu_{\ell p j}\right)=\frac{1}{(k+1)q} f\left(s, \nu_j\right),\mbox{ for all }\ell\in [k+1],p\in [q],\mbox{ and }j\in [n],
%     \label{eq:construct_f'_s_alpha}\\
%     &\hat{f}\left(\lambda_{\ell i}, \mu_{p j}\right)=\frac{1}{(k+1)q} f\left(\lambda_i, \mu_{j}\right),\mbox{ for all } \left(\lambda_{\ell i}, \mu_{p j}\right)\in \hat{\cE},
%     \label{eq:construct_f'_beta_gamma}
%     \\
%     & \hat{f}\left(\nu_{\ell p i}, \mu_{p j}\right)=\frac{1}{(k+1)q} f\left(\nu_i, \mu_{j}\right)
%     \mbox{ for all }\left(\nu_{\ell p i}, \mu_{p j}\right)\in \hat{\cE},\label{eq:construct_f'_alpha_gamma}
%     \\
%     & \hat{f}\left(\mu_{pj}, t\right)=\frac{1}{q}f\left(\mu_{j}, t\right),\mbox{ for all } p\in [q]\mbox{ and } j\in [n].
%     \label{eq:construct_f'_gamma_t}
% \end{align}
%\end{subequations}
We need the following lemma:

\begin{Lemma}
\label{lem:f_hat_flow}
    The map $\hat{f}$ defined in~\eqref{eq:construct_f_hat} is a flow.
\end{Lemma}

\begin{proof}
We follow the same procedure as in Lemma~\ref{lem:f_flow} and show that $\hat{f}$ satisfies the conservation of flow and the capacity constraint. 
\xc{Proof that $\hat{f}$ satisfies the conservation of flow.} We show that~\eqref{eq:flow_balance} holds for every node in $\hat \cV_L\cup \hat \cV_R$:
\begin{enumerate}
    \item  For any given $\lambda_{\ell i}\in \hat \cV_L$, the only incoming neighbor is $s$ and the set of outgoing neighbors is 
    $$
    N_{\rm out}(\lambda_{\ell i}) = \{(\lambda_{\ell i}, \mu_{pj}) \mid p \in [q] \mbox{ and }  (\lambda_i,\mu_j)\in \cE  \}. 
    $$ 
    It follows that
    \begin{equation}\label{eq:verify_hat_f_s_beta}
    \hat{f}\left(s, \lambda_{\ell i}\right)=\frac{1}{k} f\left(s, \lambda_i\right)
    =\frac{1}{k}\sum_{j=1}^n f\left(\lambda_i, \mu_{j}\right)
    =q\sum_{j=1}^n \hat{f}(\lambda_{\ell i}, \mu_{p j})
    =\sum_{p=1}^{q}\sum_{j=1}^n \hat{f}(\lambda_{\ell i}, \mu_{p j}),
    \end{equation}
    where the first equality follows from~\eqref{eq:def_edge_set},~\eqref{eq:def_phi}, and~\eqref{eq:f_hat_s_L}, the second equality holds due to the conservation of flow at node $\lambda_i$, for any $i\in [m]$, the third equality follows  from~\eqref{eq:def_edge_set},~\eqref{eq:def_phi}, and~\eqref{eq:f_hat_L_R}, and the last equality holds because all $\hat f(\lambda_{\ell i},\mu_{pj})$, for $\ell \in [k]$ and $p\in [q]$, are identical with each other (see~\eqref{eq:f_hat_L_R}).
    
    \item  For any given $\nu_{\ell p i}\in \hat \cV_L$, the only incoming neighbor is $s$ and the set of outgoing neighbors is
    $$
    N_{\rm out}(\nu_{\ell p i}) = \{(\nu_{\ell p i}, \mu_{pj}) \mid   (\nu_i,\mu_j)\in \cE  \}. 
    $$  
    We then have that
    \begin{align*}
    \hat{f}\left(s, \nu_{\ell p i}\right)=\frac{1}{kq}f(s,\nu_i)
    =\frac{1}{kq}\sum_{j=1}^n f(\nu_i,\mu_{j}) 
    =\sum_{j=1}^n\hat{f}\left(\nu_{\ell p i}, \mu_{p j}\right),
    %\\=\frac{1}{(k+1)q}(k+1)q\sum_{j=1}^n\hat{f}\left(\nu_{\ell p i}, \mu_{p j}\right)
    \end{align*}
where the first equality follows from~\eqref{eq:def_edge_set},~\eqref{eq:def_phi} and~\eqref{eq:f_hat_s_L}, the second equality holds due to the conservation of flow at node $\nu_i$, for any $i\in [n]$, and the last equality follows from~\eqref{eq:def_edge_set},~\eqref{eq:def_phi} and~\eqref{eq:f_hat_L_R}.

    \item  Finally, for any $\mu_{pj}\in \hat \cV_R$, the only outgoing neighbor is the sink $t$ and the set of incoming neighbors is
    \begin{align*}
    N_{\rm in}(\mu_{p j}) = & \{(\nu_{\ell p i}, \mu_{pj}) \mid \ell\in [k] \mbox{ and }  (\nu_i,\mu_j)\in \cE  \} \,  \\  &\cup
    \{(\lambda_{\ell i}, \mu_{pj}) \mid \ell\in [k] \mbox{ and } (\lambda_i,\mu_j)\in \cE  \}. 
    \end{align*}
    We obtain that
    \begin{align*}
    \hat{f}(\mu_{pj},t)
    =&
    \frac{1}{kq}\sum_{\ell=1}^{k} f(\mu_{j},t)
    \\
    =&
    \sum_{\ell=1}^{k} \sum_{i=1}^m \frac{1}{kq}f(\lambda_i,\mu_{j})
     +\sum_{\ell=1}^{k} \sum_{i=1}^n \frac{1}{kq}f(\nu_i,\mu_{j})
    \\
    =&
    { \sum_{\ell=1}^{k} \sum_{i=1}^m \hat{f}(\lambda_{\ell i}, \mu_{p j})}+\sum_{\ell=1}^{k} \sum_{i=1}^n \hat{f}\left(\nu_{\ell p i}, \mu_{p j}\right),
\end{align*}
where the first equality follows from~\eqref{eq:def_edge_set},~\eqref{eq:def_phi}and~\eqref{eq:f_hat_R_t}, the second equality holds due to the conservation of flow at node $\mu_j$, for any $j\in [n]$,  and the last equality follows from~\eqref{eq:def_edge_set},~\eqref{eq:def_phi} and~\eqref{eq:f_hat_L_R}.
\end{enumerate}

\xc{Proof that $\hat f$ satisfies the capacity constraint.}
We show below that $\hat f(e) \leq \hat c(e)$ for all $e\in \hat \cE$, with $\hat c$ defined in~\eqref{eq:def_capacity_f_bar}. First, we have that
\begin{align}
\label{eq:f_hat_s_beta}
\hat{f}\left(s, \lambda_{\ell i}\right)&=\frac{1}{k} f\left(s, \lambda_i\right)\leq \frac{1}{k}\cdot k=1,
\\
\hat{f}\left(s, \nu_{\ell p i}\right)&=\frac{1}{kq} f\left(s, \nu_i\right)\leq 1, \notag
\\
\hat{f}\left(\nu_{\ell p i}, \mu_{p j}\right)&=\frac{1}{kq} f\left(\nu_i, \mu_{j}\right)\leq 1, \notag
\\
\hat{f}\left(\mu_{pj}, t\right)&= \frac{1}{q} f\left(\mu_{j}, t\right)\leq 1,\notag
\end{align}
where the equalities follow from~\eqref{eq:def_edge_set},~\eqref{eq:def_phi}, and~\eqref{eq:construct_f_hat} and the inequalities follow from the fact that $f(e)\leq c(e)$ for all $e\in \cE$, with $c(e)$ defined given in~\eqref{eq:def_capacity_f}. 
Finally, for the edges $(\lambda_{\ell i}, \mu_{p j})\in \hat \cE$, we have that 
$$\hat{f}\left(\lambda_{\ell i}, \mu_{p j}\right)\leq \hat{f}\left(s,\lambda_{\ell i}\right)\leq 1,$$
where the first inequality follows from~\eqref{eq:verify_hat_f_s_beta} and the second inequality follows from~\eqref{eq:f_hat_s_beta}.
\end{proof}

With the above lemma, we show that $\theta(k, q) \leq \hat{\theta}(k, q)$. Let $f$ be given as in~\eqref{eq:construct_f_hat}. By Lemma~\ref{lem:f_hat_flow}, $\hat f\in \hat{ \mathcal{F}}$. Moreover, 
$$|\hat{f}|=\sum_{p=1}^{q}\sum_{j=1}^{n}\hat{f}\left(\mu_{p j}, t\right)=\sum_{j=1}^n f(\mu_j,t)=|f|,$$
where the first and the last equality follow from definition~\eqref{eq:value}, and the second equality follows from~\eqref{eq:f_hat_R_t} and~\eqref{eq:def_phi}.\unskip\nobreak\hfill\mbox{\qed}

\section{Conclusions and Outlooks}\label{sec:conclusion}
In this paper, we have addressed structural controllability for ensembles of switched, sparse linear systems where all individual systems share a common sparsity pattern $\bS$. 
For a prescribed pair \((k,q)\), with \(k\) the number of subsystems and \(q\) the number of individual systems, we derived a necessary and sufficient condition for $\bS$ to be structurally controllable. We have also provided a complete characterization of the sparsity patterns $\bS$ that are structurally controllable for $(k,q)$, for some finite $k$ and for all \(q \in \mathbb{N}\). We have further presented polynomial algorithms for checking the conditions and for computing the smallest value $k^*$. These results are given in Section~\ref{sec: Problem and Result}, and their proofs are given in Sections~\ref{sec:proof_main} and~\ref{sec:proof_flow_equiv}. {  In our future work, we will extend the results to the heterogeneous setting where the sparsity patterns for different individual systems and/or for different subsystems are not necessarily the same.}

%\newpage
% \printbibliograph

%\newpage
%\input{lemma3_apendex}

%\nocite{*}
\printbibliography
\end{document}